\documentclass[a4paper,11pt]{amsart}
\usepackage{amssymb}
 \usepackage{color}
 \usepackage{hyperref}

\newcommand{\R}{\mathbb{R}}

\newcommand{\bX}{\mathbb{X}}
\newcommand{\bS}{\mathbb{S}}
 
 \newcommand{\cB}{\mathcal{B}}
  \newcommand{\cH}{\mathcal{H}}
  \newcommand{\cI}{\mathcal{I}}
    \newcommand{\cJ}{\mathcal{J}}
 \newcommand{\cN}{\mathcal{N}}
 \newcommand{\vs}{\smallskip}

\newtheorem{theorem}{Theorem}[section]

\newtheorem{proposition}[theorem]{Proposition}

\newtheorem{lemma}[theorem]{Lemma}
\newtheorem{corollary}[theorem]{Corollary}
\newtheorem{remark}[theorem]{Remark}

 \usepackage{cite}

\begin{document}

\begin{title}[Improved Leray inequality]{\bf  Improved Leray inequality and Trudinger--Moser Type\\[1mm]
 Inequalities  involving the Leray potential   }
\end{title}
\begin{author}[H. Chen, Y. Du and F. Zhou]{}
\end{author}
\date{\today}
\maketitle


\begin{center}
 {\small
 Huyuan Chen \footnote{Email: chenhuyuan@yeah.net} 

\medskip
Department of Mathematics, Jiangxi Normal University,\\
Nanchang, Jiangxi 330022, PR China \\[14pt]

Yihong Du\footnote{Email: ydu@une.edu.au}\medskip

 School of Science and Technology, University of New England, \\ Armidale, NSW 2351, Australia \\[14pt]
 
 Feng Zhou\footnote{Email: fzhou@math.ecnu.edu.cn} 
 
 \medskip

{\small Center for PDEs, School of Mathematical Sciences, East China Normal University,\\
Shanghai Key Laboratory of PMMP, Shanghai 200062, PR China }  \\[16pt]
}  
 
\begin{abstract}
We obtain three types of results in this paper. Firstly we improve Leray's inequality by providing several types of reminder terms, secondly we introduce several Hilbert spaces based on these improved Leray inequalities and discuss their embedding properties, thirdly  we 
  obtain some  Trudinger-Moser type inequalities in the unit ball of $\R^2$ associated with the norms of these Hilbert spaces where the Leray potential is used.
   Our approach is based on analysis of radially symmetric functions.

 \vspace{4mm}
  \noindent {\small {\bf Keywords}:    Leray Inequality; Trudinger-Moser Inequalities; Leray Potential. }\vspace{1.5mm}

\noindent {\small {\bf MSC2010}: 	46E35; 26D15. }

\end{abstract}
\end{center}


\vspace{3mm}

\setcounter{equation}{0}
\section{Introduction}

The well-known Leray inequality states: 
\begin{equation}\label{ei 1.1}
\int_{B_1}|\nabla w|^2 dx -\frac14 \int_{B_1}\frac{w^2}{|x|^2(\ln\frac1{|x|})^2}dx>0,\quad\forall\, w\in C^\infty_c(B_1),\, w\not\equiv 0,
\end{equation}
where $B_1$ is the unit ball in $\R^2$, and more generally, we will write $B_r=B_r(0)$ with $B_r(x_0)$ denoting the open ball with radius $r$ centered at $x_0$ in $\mathbb R^2$. Leray used this inequality  in his study
of two dimensional viscous flows \cite{L}.

Thanks to the logarithmic function, the potential function  $\frac{1}{|x|^2(-\ln |x|)^2}$ in \eqref{ei 1.1}
has a weaker singularity at the origin than the usual inverse-square potential $1/|x|^2$. On the other hand, since $-\ln |x|\sim 1-|x|$ near $|x|=1$,
the potential function $\frac{1}{|x|^2(-\ln |x|)^2}$ has a singularity of order $(1-|x|)^{-2}$ at the boundary $\bS^1=\partial B_1$. Note that $\frac14$ is a critical coefficient both for   $\frac{1}{|x|^2(-\ln |x|)^2}$ at the origin and for  $\frac{1}{(1-|x|)^2}$ at the boundary. \smallskip

For fixed $R\in (0,1)$, since any function $w\in C^\infty_c(B_R)$ can be regarded as a function in $C^\infty_c(B_1)$ after extension by 0 to $B_1\setminus B_R$, it follows from
\eqref{ei 1.1} that
\begin{equation}\label{ei 1.1R}
\int_{B_R}|\nabla w|^2 dx -\frac14 \int_{B_R}\frac{w^2}{|x|^2(\ln\frac1{|x|})^2}dx>0,\quad\forall\, w\in C^\infty_c(B_R),\, w\not=0.
\end{equation}
We note that in \eqref{ei 1.1R}, the potential function in $B_R$ has only one singularity at $x=0$. For this inequality, there is an improved version with a remainder term, which follows from a more general result of Barbatis, Filippas and Tertikas  \cite{BFT} (whose left side can be reduced to the form of \eqref{ei 1.1R} after a change of variable; see\eqref{r_0-improved}):
\begin{equation}\label{R-improved}
\int_{B_1}|\nabla w|^2 dx -\frac14 \int_{B_1}\frac{w^2}{|x|^2(\ln\frac e{|x|})^2}dx\geq \frac 14 \sum_{i=2}^\infty\int_{B_1}\frac{|w|^2}{|x|^2}\prod_{j=1}^iX_j^2(|x|)dx,\quad\forall\, w\in C^\infty_c(B_1),
\end{equation}
where 
\begin{equation}\label{X_i}
\mbox{$X_1(r)=(\ln \frac er)^{-1}, \ X_k(r)=X_1(X_{k-1}(r))$ for $k=2,...$. }
\end{equation}

Under  the change of variable $x\to ex$, \eqref{R-improved} is reduced to the following equivalent form, which is consistent to \eqref{ei 1.1R} with $R= r_0:=e^{-1}$:
\begin{equation}\label{r_0-improved} \int_{B_{r_0}}\!\!|\nabla w|^2 dx -\frac14 \int_{B_{r_0}}\frac{w^2}{|x|^2(\ln\frac 1{|x|})^2}dx\geq \frac 14 \sum_{i=2}^\infty\int_{B_{r_0}}\frac{|w|^2}{|x|^2}\prod_{j=1}^i X_j^2(e|x|)dx, \ \forall\, w\in C^\infty_c(B_{r_0}).
\end{equation}

Let us note that in \eqref{ei 1.1}, the ball $B_1$ can be replaced by any domain $\Omega\subset B_1$ containing the origin, since any function $w\in C_c^\infty(\Omega)$ can be regarded as a function in $C^\infty_c(B_1)$ after extension by 0 over $B_1\setminus \Omega$. Similarly in  \eqref{r_0-improved}, the ball $B_{r_0}$ can be replaced by any domain $\Omega\subset B_{r_0}$ containing the origin.

\subsection{Leray's inequality with a remainder term}
In the first part of this paper, we obtain an improved version of \eqref{ei 1.1}, where the potential function has a double singularity (at both $x=0$ and $|x|=1$), and also an improved version of \eqref{ei 1.1R} which is different from  \eqref{r_0-improved}. More precisely we prove the following theorem.

\begin{theorem}\label{teo 3.1-Leray} The following inequalities hold:
\begin{itemize}
\item[(i)] {\rm (Leray's inequality with a remainder term)} There exists  $\mu_2>0$   such that for  any  $u\in C^1_c(B_1)$ 
\begin{equation}\label{e 3.1-2}
  \int_{B_1}\Big(|\nabla u|^2 -\frac14\frac{|u|^2}{|x|^2(\ln\frac 1{|x|})^2}\Big) dx\geq  \mu_2  \int_{B_1}  \frac{|u|^2 }{|x|^{2}  (\ln \frac 1{|x|})^2 \big(1+|\ln \ln \frac{1}{|x|}|\big)^2}\, dx .
\end{equation}

 \item[(ii)]  {\rm (Leray's inequality with a remainder term for radial functions)} For any $q > 2$,  there exists  $\mu_{q}>0$   such that for every  $u\in C^1_{{\rm rad}, c}(B_1)$, 
\begin{equation}\label{e 3.1-q ral}
  \int_{B_1}\Big(|\nabla u|^2 -\frac14\frac{|u|^2}{|x|^2(\ln\frac 1{|x|})^2}\Big) dx\geq  \mu_{q} \Big(\int_{B_1}  \frac{|u|^{q} }{|x|^{2}  \Big[(\ln\frac 1{ |x|})
   \big(1+|\ln \ln \frac{1}{|x|}|\big)\Big]^{1+\frac{q}{2}}}\, dx\Big)^{\frac2{q}}.
\end{equation}

 \item[(iii)] {\rm (Leray's inequality with a remainder term and singularity at 0 only)} For any $q > 2$ and $r_0=e^{-1}$,  there exists $\mu_{q}>0$  such that for  every  
 $u\in C^1_{c}(B_{r_0})$, 
\begin{equation}\label{e 3.1-q}
  \int_{B_{r_0}}\Big(|\nabla u|^2 -\frac14\frac{|u|^2}{|x|^2(\ln\frac 1{|x|})^2}\Big) dx\geq  \mu_{q} \Big(\int_{B_{r_0}}  \frac{|u|^{q} }{|x|^{2}  \Big[\big(\ln \frac 1{|x|}\big)
   \big(1+|\ln \ln \frac{1}{|x|}|\big)\Big]^{1+\frac{q}{2}}}\, dx\Big)^{\frac2{q}}.
\end{equation}
\end{itemize}
\end{theorem}

Motivated by the Leray inequality (\ref{ei 1.1}),  we denote by $\cH^1_{\mu, 0}(B_1)$, for $\mu\geq -\frac14$,   the completion of $C_c^{\infty}(B_1)$ under the norm
$$
\|u\|_\mu=\sqrt{\int_{B_1}  \Big( |\nabla u |^2dx  +\mu  \frac{u^2}{|x|^2(-\ln|x|)^2}\Big) dx}, 
$$
and so $\cH^1_{\mu, 0}(B_1)$  is a Hilbert space with inner product
$$
\langle u,v\rangle_\mu=\int_{B_1}  \Big(  \nabla u \cdot\nabla v\,dx   +\mu  \frac{u v }{|x|^2(-\ln|x|)^2} \Big)  dx. 
$$
Set $$
\cH_0^1(B_1)=\cH^1_{0, 0}(B_1) \quad{\rm and}\quad \hat \cH_0^1(B_1)=\cH^1_{-\frac14, 0}(B_1).
$$ 
For $q\in [1,+\infty)$, as usual we denote by $W^{1,q}_0(B_1)$   the completion of $C_c^\infty(B_1)$ under  the norm
$$
\|u\|_{W^{1,q}}=\Big(\int_{B_1}   |\nabla u |^q  dx\Big)^{\frac1q},
$$
and denote $\cH_0^1(B_1)=W^{1,2}_0(B_1)$.
We will show that
 \begin{equation}\label{mu>=-1/4}\mbox{
 $\cH^1_{\mu, 0}(B_1)= \cH^1_{0}(B_1)$ for $\mu>-\frac14$,  but  $ \cH^1_{0}(B_1)\varsubsetneqq \hat\cH_{0}^1(B_1)$.	  
}\end{equation}

\smallskip

\subsection{Some embedding results}

The second part of this paper gives some embedding results for the Hilbert space $\hat\cH_0^1(B_1)$.

 \begin{theorem}\label{embedding}
$(i)$ The embedding $\hat\cH_0^1(B_1)\hookrightarrow  W^{1,q}_0(B_1)$ is continuous for any $q\in[1,\,  2)$.\vs

$(ii)$    The embedding $\hat\cH_0^1(B_1)\hookrightarrow  L^p(B_1, |x|^{-\alpha} dx)$ is compact for any $p\in[ 1,+\infty)$ and  $\alpha\in [0,2)$.
 \end{theorem}
 
 The following embeddings are compact and involve more general weight functions. 
 
  \begin{theorem}\label{embedding-t}
  Let $V:\, B_1\setminus\{0\}\to [0,+\infty)$ be a continuous nonzero potential. 
   
 \begin{itemize}
\item[(i)] If
$$\begin{cases}
\displaystyle \lim_{|x|\to0^+} V(x) |x|^{2} (-\ln |x|)^{2} \big(\ln \ln \frac{1}{|x|}\big)^{2}=0,\\
\displaystyle \lim_{|x|\to1^-} V(x) (-\ln |x|)^{2} \big(\ln \ln \frac{1}{|x|}\big)^{2}=0, 
\end{cases}
$$
then 
 the embedding $\hat\cH_0^1(B_1)\hookrightarrow  L^2(B_1,  V dx)$ is compact.
 
 \item[(ii)] If $q>2$ and
$$\lim_{|x|\to0^+} V(x) |x|^{2} \Big[(-\ln |x|) \big(\ln \ln \frac{1}{|x|}\big)\Big]^{1+\frac q2}=0, $$
  then for any $r\in (0,1)$, the embedding $\hat\cH^1_0(B_{r})\hookrightarrow  L^q(B_{r}, Vdx)$ is compact.
 \end{itemize}
 \end{theorem}

 \begin{remark}\label{rm em1}
  It follows in particular that
  \begin{itemize}
\item[(i)] for any $\beta>1$,   
 the following embedding is compact:
 \[
 \hat\cH^1_0(B_1)\hookrightarrow  L^2\big(B_1, |x|^{-2} (-\ln |x|)^{-2} \big(1+ |\ln \ln \frac{1}{|x|}|\big)^{-2\beta} dx\big).
 \] 
Moreover, the embedding inequality   \eqref{e 3.1-2} holds for $u\in \hat\cH^1_0(B_1)$.

\item[(ii)] For any $q>2$, $r\in (0,1)$ and $\beta>1$,   
 the  following embedding is compact:
 \[
 \hat\cH^1_0(B_{r})\hookrightarrow  L^q\big(B_{r}, |x|^{-2} (-\ln |x|)^{1+\frac q2} \big(1+ |\ln \ln \frac{1}{|x|}|\big)^{-\beta(1 +\frac q2)} dx\big).
 \]
 Moreover, the embedding inequality   \eqref{e 3.1-q} holds for $u\in \hat\cH^1_0(B_1)$.   
  \end{itemize}
  \end{remark}

The above embedding results and their proofs can be used to obtain the following conclusions.

 \begin{theorem}\label{coro 2.1} 
$(i)$ Let  $w\in C_c^{1}(B_1)$,    then  there exists $c>0$ independent of $w$ such that
\begin{equation}\label{eqn change1}
\int_{B_1}  |\nabla w|^2  \ln \frac 1{|x|} dx\geq c \int_{B_1} |x|^{-2} (\ln \frac 1 {|x|} )^{-1} \big(1+\big|\ln \ln \frac{1}{|x|}\big|\big)^{-2} w^2 dx 
\end{equation}
and for $r\in(0,1)$
\begin{equation}\label{eqn change1-q}
\int_{B_r}  |\nabla w|^2  \ln \frac 1{|x|} dx\geq c \int_{B_r} |x|^{-2} (\ln \frac 1 {|x|} )^{-1} \big(1+\big|\ln \ln \frac{1}{|x|}\big|\big)^{-1-\frac q2} w^q dx.
\end{equation}

$(ii)$
If $\cH^1_{\ln, 0}(B_r)$ is  the completion of $C_c^{\infty}(B_r)$ under the norm
$$
\|w\|_{\ln}=\sqrt{\int_{B_r}   |\nabla w |^2  \ln \frac 1{|x|}  dx},  
$$
and $V:\, B_1\setminus\{0\}\to [0,+\infty)$ is a continuous nonzero potential, then 

$(a)$ 
 the embedding $\cH^1_{\ln, 0}(B_1) \hookrightarrow  L^2(B_1,  V  dx)$ is compact
if 
 $$\begin{cases}
\displaystyle \lim_{|x|\to0^+} V(x) |x|^{2} (\ln \frac 1 {|x|} ) \big(1+\big|\ln \ln \frac{1}{|x|}\big|\big)^{2}=0,\\
\displaystyle \lim_{|x|\to1^-} V(x)   (\ln \frac 1 {|x|} ) \big(1+\big|\ln \ln \frac{1}{|x|}\big|\big)^{2}=0; 
\end{cases}
$$

$(b)$ for $r\in(0,1)$, 
 the embedding $\cH^1_{\ln, 0}(B_r) \hookrightarrow  L^q(B_r,  V  dx)$ is compact
if 
$$ 
\displaystyle \lim_{|x|\to0^+} V(x) |x|^{2} (\ln \frac 1 {|x|} ) \big(\big|\ln \ln \frac{1}{|x|}\big|\big)^{1+\frac q2}=0.
$$

\end{theorem}
 
 
\subsection{Trudinger-Moser  type inequalities}
  
  The Trudinger-Moser inequality (\!\!\cite{Moser, Tr})
\begin{equation}\label{TM}
\sup_{ \int_\Omega|\nabla u|^2dx\leq 1,\ u\in C_c^\infty(\Omega)}\int_\Omega e^{4\pi u^2} dx<\infty,
\end{equation}
where $\Omega$ is a bounded domain in $\R^2$, is an analogue of the following limiting Sobolev inequality in dimensions $N\geq 3$:
\[
\sup_{\int_{\R^N}|\nabla u|^2dx\leq 1,\ u\in C_c^\infty(\R^N)}\int_{\R^N}|u|^{2^*}  dx<\infty,\ 2^*=\frac{2N}{N-2}.
\]
There are several extensions of \eqref{TM} in the literature, where the constraint $\int_\Omega|\nabla u|^2dx\leq 1$ is replaced by
\[
\int_\Omega \Big[|\nabla u|^2-V(x) u^2\Big]dx\leq 1
\]
with a suitable potential function $V(x)$. For example, in Wang and Ye \cite{WY}, it was shown that \eqref{TM} remains valid for $\Omega=B_1$ if 
$\int_\Omega|\nabla u|^2dx\leq 1$ is replaced by 
\[
\int_{B_1} \Big[|\nabla u|^2-V_1 u^2\Big]dx\leq 1 \mbox{ with } V_1:=(1-|x|^2)^{-2}.
\]
The potential function $V_1$ here is related to the Leray potential 
\[
V_{\rm Leray}:=\frac 1{4(\ln |x|)^2(\ln\frac 1{|x|})^2}
\]
 in the following way:
\[
\lim_{r\to 1^-}  V_1(r)/V_{\rm Leray}(r)=1.
\]
Tintarev \cite{T} showed that \eqref{TM} remains valid if $\Omega=B_1$ and 
$\int_\Omega|\nabla u|^2dx\leq 1$ is replaced by 
\[
\int_{B_1} \Big[|\nabla u|^2-V_2 u^2\Big]dx\leq 1 \mbox{ with } V_2:=\frac{V_{\rm Leray}(|x|)}{\max\{\sqrt{-\ln {|x|}}, 1\}}.
\]
Clearly
\[
\lim_{|x|\to 1^-} V_2(|x|)/V_{\rm Leray}(|x|)=1,\ \lim_{|x|\to 0^+} V_2(|x|)/V_{\rm Leray}(|x|)=0.
\]
Let us note that $V_1$ is regular at the origin with $V_1(0)=1$ while $V_2$ is singular at the origin.
Set
\[
V_3(r):=\frac{1}{4r^2(1-\ln r)^2}.
\]
Then it is easily seen that $V_3(r)$ is decreasing in $r$ for $r\in (0, 1]$, $V_3(1)=1/4$ and
\[
\lim_{r\to 0^+} V_3(r)/V_{\rm Leray}(r)=1.
\]
For $V=V_3$, it was shown by Psaradakisa and Spectora \cite{PS} 
 that for any domain $\Omega\subset B_1$, and any $\epsilon>0$, there exists $A_\epsilon>0$ such that
\[
\sup_{ \int_\Omega[|\nabla u|^2-V_3 u^2]dx\leq 1,\ u\in C_c^\infty(\Omega)}\int_\Omega e^{A_\epsilon u^2(x)/(1-\ln|x|)^{\epsilon}} dx<\infty,
\]
and the result is false for $\epsilon=0$.
Mallick and Tontarev \cite{MT} subsequently proved that for any domain $\Omega\subset B_1$,  there exists $A>0$ such that
\[
\sup_{ \int_\Omega[|\nabla u|^2-V_3 u^2]dx\leq 1,\ u\in C_c^\infty(\Omega)}\int_\Omega e^{A u^2(x)/E(|x|)} dx<\infty,
\]
where
\[
E(|x|):=1-\ln(1-\ln |x|).
\]

In a different direction, Adimurthi and Druet \cite{A} proved the following result: For any bounded domain $\Omega\subset\R^2$,
\[
\sup_{ \int_\Omega|\nabla u|^2dx\leq 1,\ u\in C_c^\infty(\Omega)}\int_\Omega e^{4\pi u^2(1+\alpha \|u\|_2)} dx\begin{cases}<\infty \mbox{ if } \alpha\in [0,\lambda_1(\Omega)),\\[2mm]
=\infty \mbox{ if } \alpha\geq \lambda_1(\Omega),
\end{cases}
\]
where $\lambda_1(\Omega)$ stands for the first eigenvalue of $-\Delta$ in $H_0^1(\Omega)$, and $\|u\|_2=(\int_\Omega u^2dx)^{1/2}$.
\medskip

Inspired by these results, we will prove some related but different Trudinger-Moser type inequalities involving a potential $V$ behaving like the Leray potential near 0 and/or near $\partial B_1$.

Let $\Omega\subset B_1$ be a bounded  domain containing the origin,  $\mu\geq-\frac14$, and $V: (0,1)\to [0,\infty)$ a continuous function such that
\[
\mu V(r)\geq -V_{\rm Leray}(r)=-\frac 1{4r^2 (\ln\frac 1 r)^2} \mbox{ for } r\in (0,1).
\]
We denote by  $\cH^1_{_V,\mu, 0}(\Omega)$   the completion of $C_c^\infty(\Omega\setminus\{0\})$ under the norm
 \begin{eqnarray*}  \label{norm V}
 \|u\|_{_V,\mu}=\sqrt{\int_{\Omega}  \Big( |\nabla u |^2dx  +\mu Vu^2  \Big) dx},
 \end{eqnarray*} 
which is a Hilbert space with  inner product
$$\langle u,v\rangle_{_V,\mu}=\int_{\Omega}  \Big(  \nabla u \cdot\nabla v\,dx   +\mu V u v  \Big)  dx. $$
 
 For $\mu\geq -\frac14$,  we denote 
\begin{equation*}\label{optimal bound}
m_\mu :=4 \pi     \sqrt{1+4\mu} . 
\end{equation*}

  \begin{theorem}\label{cr 0}
Assume that $\mu>0$ and $V:(0,1)\to[0,+\infty)$ is a continuous  function satisfying 
\begin{equation}\label{con V-0}
 V(r)\geq \frac{1}{r^2(-\ln r)^2}\quad {\rm in}\ (0,1).
\end{equation}
Then the following conclusions hold:
 \begin{itemize}
 \item[(i)] For {\bf radially symmetric} functions in $\cH^1_{_V,\mu, 0}(B_1)$ we have
$$
\sup_{u {\rm \, is\,  radial}, \|u\|_{V,\mu}\leq 1}\int_{B_1} e^{m_\mu |u|^2}dx<\infty,
$$
and this result is optimal: If $\alpha>m_\mu$ and
\begin{equation}\label{con V-l}
 \lim_{r\to0^+}V(r) {r^2(-\ln r)^2}=1,
\end{equation}
then
 there exists a sequence of radially symmetric functions  which concentrate at the origin such that $\|u_n\|_{_V,\mu}\leq 1$ and 
$$\int_{B_1} e^{\alpha |u_n|^2}dx\to \infty\quad{\rm as}\ \, n\to+\infty.$$

\item[(ii)]  For general functions in $\cH^1_{_V,\mu, 0}(B_1)$ we have
$$
\sup_{ \|u\|_{V,\mu}\leq 1}\int_{B_1} e^{4\pi |u|^2}dx<\infty,
$$
and this result is optimal:
If   $\alpha>4\pi$ and \eqref{con V-l} holds, 
then
 there exists  a sequence of functions concentrating at some point away from the origin, such that  $\|u_n\|_{_V, \mu}\leq 1$  and
$$\int_{B_1} e^{\alpha |u_n|^2}dx\to +\infty\quad{\rm as}\ \, n\to+\infty.$$
\end{itemize}

 \end{theorem}
 
 Here a sequence $\{u_n\}$ is said to be concentrating  at some point $x_0$, if for any $r\in(0,1)$ and any $\epsilon>0$
 there exists $n_0>0$ such that 
 $$\int_{B_1\setminus B_r(x_0)}\Big( |\nabla u |^2dx  +\mu Vu^2  \Big) dx <\epsilon. $$

 Next we consider the case $\mu\in(-\frac14,0)$.

 \begin{theorem}\label{cr 1}
 Suppose $\mu\in(-\frac14,0)$ and  $V:   (0,1)\to[0,+\infty)$ is continuous and verifies
 \begin{equation}\label{con V}
 V(r)\leq \frac{1}{r^2(-\ln r)^2}.
\end{equation}
  Then the following conclusions hold:
  \begin{itemize}
  \item[(i)]
  For  {\bf radially symmetric} functions,
  $$
\sup_{u {\rm \, is\,  radial}, \|u\|_{V,\mu}\leq 1}\int_{B_1} e^{m_\mu |u|^2}dx<\infty,
$$
and this result is optimal: If $\alpha>m_\mu$ and \eqref{con V-l} holds,
then
 there exists a sequence of radially symmetric functions  which concentrate at the origin such that $\|u_n\|_{_V,\mu}\leq 1$ and 
$$\int_{B_1} e^{\alpha |u_n|^2}dx\to \infty\quad{\rm as}\ \, n\to+\infty.$$
    \item[(ii)] For general functions, if $V$ is {\bf decreasing} in $(0,1)$ and verifies \eqref{con V},
 then 
 $$
\sup_{\|u\|_{V,\mu}\leq 1}\int_{B_1} e^{m_\mu |u|^2}dx<\infty.
$$
\item[(iii)] The result in {\rm (ii)} is optimal:  If 
 \eqref{con V-l} holds,   
then  for any $\alpha>m_\mu$, there exists a sequence $\{u_n\}_n$ concentrating at the origin such that  $\|u_n\|_{_V,\mu}\leq 1$ and 
$$\int_{B_1} e^{\alpha |u_n|^2}dx\to \infty\quad{\rm as}\ \, n\to\infty.$$
\end{itemize}
 \end{theorem}
 
 Let us note that $V(r):=\frac{1}{r^2(1-\ln r)^2}$ is decreasing in $(0,1]$ and satisfies both \eqref{con V-l} and \eqref{con V}.

    \begin{corollary}\label{cr 2-real}
  Let    $\mu\in(-\frac14,0)$, $r_0=1/e$ and $V(r)=\frac{1}{r^2(\ln \frac 1r)^2}$. Then
$$\sup_{\|u\|_{\cH^1_{_V,\mu, 0}(B_{r_0})}\leq 1}\int_{B_{r_0}} e^{m_\mu |u|^2}dx<\infty ,$$
and the exponent $m_\mu$ is optimal.
 \end{corollary}
 \begin{proof}
 Under the change of variable $x\to ex$, the inequality is changed to an equivalent  one over $B_1$ with $V(r)$ replaced by $V(r/e)$, which is decreasing over $(0,1)$.
 Therefore we can use Theorem \ref{cr 1} (ii) to conclude.
 \end{proof}

Finally we consider the critical case  $\mu=-\frac14$.

  \begin{theorem}\label{cr 2}
  Suppose that $\mu=-\frac14$ and $V\in C((0,1))$ is nonnegative  and verifies \eqref{con V}.
Then the following conclusions hold:
\begin{itemize}
\item[(i)] For   {\bf radially symmetric} functions and
$p\in(0,1)$, $\alpha>0,$
$$\sup_{u {\rm\, is\; radial},\; \|u\|_{V,-1/4}\leq 1}\int_{B_1} e^{\alpha |u|^p}dx<\infty.$$
\item[(ii)] For general functions, 
if $V$ is {\bf decreasing} in $(0,1)$,    
then for
$p\in(0,1)$ and $\alpha>0$,
$$\sup_{\|u\|_{V,-1/4}\leq 1}\int_{B_1} e^{\alpha |u|^p}dx<\infty.$$
 \item[(iii)] If  there exist $\theta>0$ and $C>0$ such that  
\begin{equation}\label{con V-l2}
|V(r) {r^2(-\ln r)^2}-1|\leq C(-\ln r)^{-\theta} \quad \mbox{for}\ r\in(0,\frac14),
\end{equation}
then 
   there exists a sequence $\{u_n\}\subset \cH^1_{V, -1/4,0}(B_1)$ such that  $\|u_n\|_{\cH^1_{V,-1/4}(B_1)}= 1$ and
   for any $p\geq1$ and any $\alpha>0$,
$$\int_{B_1} e^{\alpha |u_n|^p }dx\to \infty\quad{\rm as}\ \, n\to\infty.$$
\end{itemize}
 \end{theorem}
 
 \begin{remark} We end this subsection with some remarks:
 \begin{itemize}
 \item[(i)] Theorem \ref{cr 0} {\rm (ii)} contrasts sharply with Theorem \ref{cr 0} {\rm (i)} and Theorem \ref{cr 1} {\rm (ii)} in that, the best components in the latter results are both given by $m_\mu$ while that in Theorem \ref{cr 0} {\rm (ii)} is $4\pi \;(=m_0)$. 
  \item[(ii)] By Corollaries \ref{cr:main-0} and \ref{cr:main-1} below,  under  condition \eqref{con V}, for $\mu>-1/4$, $\cH^1_{V,\mu, 0}(B_1)=\cH^1_0(B_1)$, but
 $ \cH^1_{0}(B_1)\varsubsetneqq \cH^1_{_V,-1/4, 0}(B_1)$ if additionally \eqref{con V-l2} holds. These clearly imply \eqref{mu>=-1/4}.

  \end{itemize}
 
 \end{remark}
 
\subsection{Organisation of the paper}

 The rest of the paper is organised as follows. In Section 2, 
we prove Theorem \ref{teo 3.1-Leray} based on a general inequality \cite[Theorem 3, p.44]{M}. In Section 3, we consider the embedding properties of $\hat\cH^1_0(B_1)$ and related spaces,
where we in particular prove Theorems \ref{embedding},  \ref{embedding-t} and \ref{coro 2.1}. Section 4 is devoted to the proof of Trudinger--Moser type inequalities for radially symmetric functions, where rather subtle and involved analysis is used to obtain the desired results; this  is perhaps the most technically demanding part of the paper. Section 5 extends the results of Section 4 for radial functions to general functions, where significant difference is  revealed (see Remark 1.8 (i) above). With all the main ingredients ready, the proofs of Theorems \ref{cr 0}, \ref{cr 1} and \ref{cr 2} are completed at the end of Section 5.
 
 \setcounter{equation}{0}
\section{Leray's inequality with remainder terms  }

We prove Theorem \ref{teo 3.1-Leray} in this section. 
Our analysis will be based on the following proposition. 

\begin{proposition}\label{teo 3.1-0}{\rm \cite[Theorem 3, p.44]{M}}
Suppose $1\leq p\leq q\leq \infty$, $\gamma$ and $\nu$ are measures such that 
$$\cB=\sup_{r>0} \gamma\big((0,r)\big)^{\frac1q} \Big(\int_r^\infty (\frac{d\nu^*}{dr})^{-\frac{1}{p-1}}\Big)^{\frac{p-1}{p}}<+\infty,$$
where    $\nu^*$ is the absolutely continuous part of $\nu$. 
 Then there exists $c>0$ such that for any $f\in C(\R_+)\cap  L^p(\R_+,d\nu)$  
$$\Big(\int_0^\infty \Big|\int_0^r f(t)dt\Big|^q d\gamma(r)\Big)^{\frac1q}
\leq c \Big(\int_0^\infty |f(r)|^pd\nu(r)\Big)^{\frac1p}.$$
Moreover, if $c$ is the best constant in the above inequality, then 
$$\cB\leq c\leq \cB (\frac{q}{q-1})^{\frac{p-1}{p}} q^{\frac1q},$$
and $c=\cB$ if $p=1$ or $q=\infty$.
\end{proposition}

\begin{lemma}\label{lm 3.1} Suppose that $q\geq 2$, $a_0:=e^{-e}$,  and the functions $a(r)$ and $b(r)$ satisfy
\[
 a(r)\geq \beta r(-\ln r),\ \  b(r)\geq \beta r(\ln \frac1r) (\ln\ln \frac 1r)^{\frac q2+1}
\]
for all $r\in (0,a_0]$ and some constant $\beta>0$.
Then there exists $c_1=c_1(q, \beta)>0$ such that for any $v\in C_0^1\big([0, 1)\big)$,
    \begin{equation}\label{es 3.1}
\int_0^{a_0} v'(r)^2 a(r) dr\geq  c_1\Big(\int_0^{a_0}\frac{|v(r)|^q}{b(r)} dr\Big)^{\frac2q}.
    \end{equation}
\end{lemma}
\noindent{\bf Proof. }   In what follows, for any interval $I\subset\mathbb R$,  $\chi_I$ will denote the characteristic function of $I:$
   \[
   \chi_I(r)=\begin{cases} 1, & r\in I,\\
0,& r\not\in I. \end{cases}
\]
Set  $I=[0, a_0]$ and
$$d\gamma=\frac{\chi_{_I}(s)}{b(s)} ds, \quad
d\nu = \frac {s(-\ln s)} {\chi_{_I}(s)} ds. $$
Then for $r\in I$,
    \begin{eqnarray*}
 \gamma([0,r])\leq \frac 1{\beta_1} \int_0^{r}\frac{1}{s\ln\frac1s (\ln\ln\frac1s)^{\frac q2+1}} ds =  \frac2{\beta_1q} \, (\ln\ln\frac1r)^{-\frac q2} 
   \end{eqnarray*}
   and for $r>a_0$, 
   $$\gamma([0,r])=\gamma([0,a_0]) \leq \frac2{\beta_1q}.$$
   
   Moreover,  for $r\in I$,
  \begin{eqnarray*}
 \int_{r}^{+\infty} \Big(\frac{d\nu^*}{ds}\Big)^{-1}ds =  \int_{r}^{a_0}\frac1{s(-\ln s)}  ds= \ln\ln\frac1r -1
   \end{eqnarray*}
    and for $r>a_0$, 
   $$ \int_{r}^{+\infty} \Big(\frac{d\nu^*}{ds}\Big)^{-1}ds=0.$$

So we have that 
  \begin{eqnarray*}
\cB&=&\sup_{r>0} \, \gamma([0,r])^{\frac1q} \left(\int_{r}^{+\infty} \Big(\frac{d\nu^*}{ds}\Big)^{-1}ds\right)^{\frac12} 
\\[2mm]&=&\Big(\frac2{\beta_1q}\Big)^{\frac1q} \sup_{r\in I}\sqrt{\frac{\ln\ln\frac1r -1}{\ln\ln\frac1r} }
\\[2mm]&\leq&  \Big(\frac2{\beta_1q}\Big)^{\frac1q}.\  
   \end{eqnarray*}
Now we can apply Proposition \ref{teo 3.1-0} with $q\geq 2=p$ and $f(r)=\chi_I(r)v'(r)$ to obtain (\ref{es 3.1}). \hfill$\Box$\medskip

\begin{lemma}\label{lm 3.2}
Suppose $q\geq 2$ and the functions $\tilde a(r)$ and $\tilde b(r)$ satisfy
\[
 \tilde a(r)\geq \tilde \beta (-\ln r),\ \  \tilde b(r)\geq \tilde \beta (1-r) (\ln\frac 1{1-r})^{\frac q2+1}
\]
for all $r\in [a_0, 1)$ and some constant $\tilde \beta>0$. Then there exists $c_2=c_2(q,\tilde\beta)>0$ such that for any $v\in C_0^1\big((0,1)\big)$,
    \begin{eqnarray}\label{es 3.2}
\int_{a_0}^{1} |v'(r)|^2 \tilde a(r)  dr\geq  c_2\Big(\int_{a_0}^1\frac{|v(r)|^q}{\tilde b(r)} dr\Big)^{\frac2q}.
    \end{eqnarray}
\end{lemma}
\noindent{\bf Proof. } In order to use Proposition \ref{teo 3.1-0}, we set 
$$d\gamma=\frac{\chi_{_J}(s)}{\tilde b(1-s)} dr\quad {\rm and}\quad
d\nu= \frac{s}{\chi_{_J}(r)}   dr, $$
where $J=(0, 1-a_0)$. 
Then for $r\in J$,
    \begin{eqnarray*}
 \gamma([0,r])\leq \frac1{\tilde\beta} \int_{0}^{r}\frac{1}{s\big(-\ln s\big)^{\frac q2+1}} ds =  \frac2{\tilde\beta q} \,  \big(-\ln r\big)^{-\frac q2} ,
   \end{eqnarray*}
  and  for  $r\geq 1-a_0$
       \begin{eqnarray*}
 \gamma([0,r])= \gamma([0,1-a_0])\leq  \frac2{\tilde\beta q} \,  \Big(-\ln (1-a_0)\Big)^{-\frac q2} .
   \end{eqnarray*}
    
  On the other hand,
 for $r\in J$
  \begin{eqnarray*}
 \int_{r}^{+\infty} \Big(\frac{d\nu^*}{ds}\Big)^{-1}ds = \int_{r}^{1-a_0}\frac 1{s} ds= -\ln r +\ln (1-a_0)\leq -\ln r,
   \end{eqnarray*}
   and for $r\geq 1-a_0$,  
   \begin{eqnarray*}
 \int_{r}^{+\infty} \Big(\frac{d\nu^*}{ds}\Big)^{-1}ds =0.
   \end{eqnarray*}     
So we have that 
  \begin{eqnarray*}
\cB = \sup_{r>0} \, \gamma([0,r])^{\frac1q} \left(\int_{r}^{+\infty} \Big(\frac{d\nu^*}{ds}\Big)^{-1}ds\right)^{\frac12}   \leq   \Big(\frac2{\tilde\beta_1q}\Big)^{\frac1q}.\  
   \end{eqnarray*}
Now we can apply Proposition \ref{teo 3.1-0} with $q\geq 2=p$ and $f(r)=\chi_J(r)w'(r)$, $w(r)=v(1-r)$,  to obtain some $\tilde c_1>0$,
\begin{eqnarray*} 
\int_0^{1-{a_0}}|\omega'(s)|^2 s  ds\geq  \tilde c_1\Big(\int_0^{1-{a_0}} \frac{|\omega(s)|^q}{\tilde b(1-s)} ds\Big)^{\frac2q},
    \end{eqnarray*}
    i.e.
    \[
    \int^1_{{a_0}} |v'(s)|^2 (1-s)  ds\geq  \tilde c_1 \Big(\int^1_{{a_0}}  \frac{|v(s)|^q}{\tilde b(s)} ds\Big)^{\frac2q}.
    \]
Since $-\ln s\sim 1-s$ as $s\to 1^{-}$, there exists $\tilde c_2>0$ such that
\begin{eqnarray*} 
   \int^1_{{a_0}} |v'(s)|^2 \tilde a(s)  ds&\geq& \tilde\beta_1 \int^1_{{a_0}} |v'(s)|^2 (-\ln s)  ds 
   \geq \tilde c_2 \int^1_{{a_0}} |v'(s)|^2 (1-s)  ds\\
 &\geq&  \tilde c_2\tilde c_1 \Big(\int^1_{{a_0}}  \frac{|v(s)|^q}{\tilde b(s)} ds\Big)^{\frac2q}.
    \end{eqnarray*}
    Thus (\ref{es 3.2}) holds.  \hfill$\Box$\medskip\vs
 
\noindent {\bf Proof of Theorem \ref{teo 3.1-Leray}. } 
Let 
$$\cI(w)=\int_{B_1}|\nabla w|^2 dx  -\frac14 \int_{B_1}\frac{w^2}{|x|^2(-\ln|x|)^2}dx,$$
which is nonnegative by  Leray's inequality (\ref{ei 1.1}).\smallskip

{\it Part 1. }  Let $v_0$ be a radially symmetric function in $C^1_0({B_1})$, $\alpha\geq 0$ a constant,  and   
$$w_0(r)=(\alpha-\ln r)^{-\frac12} v_0(r)\quad{\rm for}\ \, r\in(0,1).$$
Then
 $w_0(0)= w_0(1)=0$, and
 we have 
\begin{eqnarray} 
  \cI_\alpha(v_0): &=& 2\pi \int^1_0 \Big(v_0'(r)^2-\frac14\frac{v_0(r)^2}{r^2(\alpha-\ln r)^2}\Big)  rdr
 \nonumber
 \\[2mm]&=&   2\pi \int^1_0 \Big(w_0'(r)^2- \frac{w_0'(r)w_0(r) }{r (\alpha-\ln r) }\Big)  r(\alpha-\ln r) dr 
 \nonumber
 \\[2mm] &=&2\pi \int^1_0  w_0'(r)^2 r(\alpha-\ln r) dr -2\pi \int^1_0 w_0'(r)w_0(r)dr\label{equivalence 1}
 \\[2mm]&=&  2\pi \int^1_0  w_0'(r)^2 r(\alpha-\ln r) dr + \pi(w_0^2(0)-w_0(1)^2)
   \nonumber\\[2mm]&=&  2\pi \int^1_0  w_0'(r)^2 r(\alpha-\ln r) dr. \nonumber
    \end{eqnarray}
 From Lemma \ref{lm 3.1} we obtain that for $q\geq 2$,
 \begin{eqnarray*} 
\int^{a_0}_0   |w_0|^2  r(-\ln r) dr &\geq&  c_1\Big(\int_0^{a_0}\frac{|w_0|^q}{r(\ln \frac1r) (\ln\ln \frac1r)^{\frac q2+1}} dr\Big)^{\frac2q},
    \end{eqnarray*}
    and by Lemma \ref{lm 3.2}
  \begin{eqnarray*} 
\int^1_{a_0}   |w_0'|^2  r(-\ln r) dr &\geq&  a_0c_2 \Big(\int_{a_0}^1\frac{\big|w_0\big|^q}{(1-r)\big(\ln \frac1{1-r}\big)^{\frac q2+1}} dr\Big)^{\frac2q}. 
    \end{eqnarray*}
    As a consequence, 
   \begin{eqnarray*}    
   \cI(v_0)=\cI_0(v_0) &=& 2\pi \int^1_0  w_0'(r)^2 r(-\ln r) dr \\[2mm] 
 &\geq&c_3 \Big( \int_0^1 \Big[\frac{\chi_{[0,a_0]}(r)}{  r (\ln \frac 1r)  (\ln \ln \frac{1}r)^{\frac {q}2+1}} +\frac{\chi_{[a_0,1]}(r)}{  (1-r)(\ln\frac{1}{1-r})^{\frac{q}2+1} } \Big]  |w_0|^qdr\Big)^{\frac2q}.
     \end{eqnarray*}
Since 
$$
\lim_{r\to1^-}\frac{r\ln \frac{1}{r}}{1-r}=\lim_{r\to1^-}\frac{-\ln\ln \frac 1r}{\ln \frac 1{1-r}}=1, \lim_{r\to 0^+} \frac{1+\ln\ln \frac 1r}{\ln\ln \frac 1r}= \lim_{r\to 1^-} \frac{1-\ln\ln \frac 1r}{-\ln\ln \frac 1r}=1,$$
we further obtain
 \begin{equation} \label{3.2-1}
 \begin{aligned}
   \cI(v_0)&\geq  c_4 \Big( \int_0^1 \Big[\frac{\chi_{[0,a_0]}(r)}{  r (\ln \frac 1r)  (1+\ln \ln \frac{1}r)^{\frac {q}2+1}} +\frac{\chi_{[a_0,1]}(r)}{  r( \ln \frac 1r) \big(1+ |\ln \ln \frac{1}r
    |\big)^{\frac{q}2+1} } \Big]  |w_0|^qdr\Big)^{\frac2q}\\
   &= c_4\Big( \int_0^1 \Big[\frac{|w_0|^q}{  r (\ln \frac 1r )(1+ |\ln \ln \frac{1}r|)^{\frac {q}2+1}}dr\Big)^{\frac2q}\\
   &=c_4\Big(\int_{B_1}  |x|^{-2}\Big[\big(\ln\frac1{|x|} \big)\big(1+\big|\ln \ln \frac{1}{|x|}\big|\big)\Big]^{-1-\frac{q}{2}}|v_0|^{q} dx\Big)^{\frac2{q}}.
      \end{aligned}
      \end{equation}
     {\it Part 2.   General functions for $q=2$ in $B_1$. }   Motivated by \cite{VZ,F}, we prove the desired inequality for $q=2$ by making use of the spherical harmonic functions. 
     
        Any $u\in C_c^\infty(B_1)$  has a decomposition
     into spherical harmonics in the following form:
     $$u(x)=\sum^\infty_{m=0}u_m(r)h_m(\sigma),$$
where  $\{h_m\}$ consists of the orthonormal eigenfunctions of the Laplacian-Beltrami operator
on the sphere  with corresponding eigenvalues $\lambda_m=m^2$ for  integer $m\geq 0$.  In particular, $u_0$ is the radial part of $u$ and $h_0=1$. 

Note that 
$$\int_{B_1} |\nabla u|^2dx= \sum^\infty_{m=0} \int_{B_1}\Big( |\nabla u_m |^2+\lambda_m\frac{u_m ^2}{|x|^2}\Big) dx $$
and  by  (\ref{3.2-1})
  \begin{eqnarray} 
  \cI(u)&=&\sum_{i=0}^{\infty}\cI(u_m) +2\pi \sum^\infty_{m=1} \Big( \lambda_m \int_0^1 \frac{ u_m^2}{|x|^2} dx\Big) \nonumber
\\[2mm] &\geq& c_7 \sum^\infty_{m=0}\Big(\int_{B_1}  |x|^{-2}\big(\ln\frac1{|x|}\big)^{-1-\frac{q}{2}}(1+ \big|\ln \ln \frac{1}{|x|}\big|)^{-q}u_m^{q} dx\Big)^{\frac2{q}} +  2\pi \sum^\infty_{m=1} \Big( \lambda_m \int_0^1\frac{ u_m^2}{r} dr\Big). \label{3.1==1}
\end{eqnarray} 
Note that   
 \begin{eqnarray*} 
 &&\int_{B_1} |x|^{-2}\big(\ln\frac1{|x|}\big)^{-2} \big(1+\big|\ln \ln \frac{1}{|x|}\big|\big)^{-2} |u|^2dx 
 \\[2mm]&=&   \int_0^1  r^{-2}\big(\ln\frac1{r}\big)^{-2} \big(1+\big|\ln \ln \frac{1}{r}\big|\big)^{-2}\int_{\bS^1} \Big|\sum^\infty_{m=0}   u_m(r)  h_m(\sigma)\Big|^2d\sigma dr
 \\[2mm]&=&   \int_0^1  r^{-2}\big(\ln\frac1{r}\big)^{-2} \big(1+\big|\ln \ln \frac{1}{r}\big|\big)^{-2} \Big|\sum^\infty_{m=0,n=0}  u_m(r) u_n (r)  \int_{\bS^1} h_m(\sigma)h_n(\sigma) d\sigma dr
 \\[2mm]&=&\sum^\infty_{m=0}  \int_{B_1}  |x|^{-2}\big(\ln\frac1{|x|}\big)^{-2} \big(1+\big|\ln \ln \frac{1}{|x|}\big|\big)^{-2} |u_m|^2 dx,
 \end{eqnarray*} 
 where we used the facts that 
 $\int_{\bS^1}h_i(\sigma)h_j(\sigma) d\sigma=0$ for $i\not=j$
 and
 $\int_{\bS^1}|h_i(\sigma)|^2 d\sigma=1$ for $i=0,1,\cdots$.
 
Now we apply (\ref{3.1==1}) with $q=2$ and obtain 
 \begin{eqnarray*} 
  \cI(u) &\geq& c_7   \sum_{i=0}^{\infty} \int_{B_1}  |x|^{-2}\big(\ln\frac1{|x|}\big)^{-2} \big(1+\big|\ln \ln \frac{1}{|x|}\big|\big)^{-2}u_m^{2} dx \nonumber
  \\[2mm] &=&c_7 \int_{B_1} |x|^{-2}\big(\ln\frac1{|x|}\big)^{-2}\big(1+ \big|\ln \ln \frac{1}{|x|}\big|\big)^{-2} u^2dx.
\end{eqnarray*} 

However,  this method  fails for $q>2$. \smallskip
     
  {\it Part 3.   General functions for $q>2$  in $B_{r_0}$. }     Let $r_q\in(0,1)$ be chosen such that the function $r^{-2}  \big[(-\ln r) \big(1+|\ln \ln \frac{1}{r}|\big)\big]^{-1-\frac{q}{2}}$ is decreasing in $(0,r_q)$.  We will make use of the  rearrangement argument.

Let $u^*$ denote the symmetric decreasing rearrangement of $u$ (extended by 0 over $\mathbb R^2\setminus B_{r_0}$), and let us denote
    $$V_1(r):=\frac{1}{r^{2}(-\ln r)^2},\quad V_2(r):=\frac{1}{r^{2}\Big[(-\ln r)\big(1+|\ln\ln\frac{1}{r}|\big)\big]^{1+\frac{q}{2}}},$$
which are decreasing in $(0,r_0)$ and in $(0, r_q)$, respectively.  By the P\'olya-Szeg\H{o} inequality we have   
$$\int_{B_{r_0}} |\nabla u^*|^2 dx\leq \int_{B_{r_0}} |\nabla u |^2 dx.
$$ 
  Combining this with (\ref{3.2-1}), using also the Hardy-Littlewood inequality and $V_1^*=V_1$ in $B_{r_0}$,  we obtain
    \begin{eqnarray*}  
   \cI(u)&=& \int_{B_{r_0}} |\nabla u |^2 dx-\frac14 \int_{B_{r_0}}  u^2 |x|^{-2}  \big(\ln\frac1{|x|}\big)^{-2} dx
    \\[2mm]&\geq& \int_{B_{r_0}} |\nabla u^* |^2 dx-\frac14\int_{B_{r_0}}  u^2  V_1   dx
    \\[2mm]&\geq & \int_{B_{r_0}} |\nabla u^* |^2 dx-\frac14\int_{B_{r_0}}  (u^*)^2  V_1   dx\\
    &\geq &   c  \Big(\int_{B_{r_0}}   (u^*)^{q} V_2dx\Big)^{\frac2{q}}.
      \end{eqnarray*}  
      
     If $r_q\geq r_0$, then
     \[
    \int_{B_{r_0}}   (u^*)^{q} V_2dx= \int_{B_{r_0}}   (u^*)^{q} V_2^*dx\geq \int_{B_{r_0}}   u^{q} V_2dx,
 \]
    which ends the proof. 
 
 If $r_q<r_0$,
then due to $V_2(0^+)=+\infty$ there exists $\hat r_0\in (0, r_q)$ such that $V_2(\hat r_0)\geq \max_{r\in [r_q, R]}V_2(r)$.
We now define
\[
\bar V_2(r):=\begin{cases} V_2(r), &r\in [0, \hat r_0),\\
V_2(\hat r_0),& r\geq \hat r_0,\end{cases} \ \  \underline {V_2}(r):=\begin{cases} V_2(r), & r\in [0, \hat r_0),\\
\min_{r\in [r_q, R]}V_2(r),& r\geq \hat r_0.
\end{cases}
\]
Then it is clear that
\[\begin{cases}
\underline {V_2}\leq V_2\leq \bar V_2,\\[1.5mm]
\underline {V_2}\geq c_q \bar V_2\ \mbox{ for some constant } c_q\in (0,1),\\[1.5mm]
\mbox{both $\underline {V_2}$ and $ \bar V_2$ are non-increasing.}
\end{cases}
\]
Now
$$ \int_{B_{R}} (u^*)^q   V_2   dx \geq \int_{B_{R}} (u^*)^q   \underline {V_2}  dx= \int_{B_{R}} (u^*)^q   \underline {V_2}^*  dx \geq   \int_{B_{R}}  u^q \underline {V_2}dx\geq  c_q\int_{B_{R}}  u^q {V_2}dx,
$$      
and thus 
 \begin{eqnarray*}  
   \cI(u)&\geq &   c  \Big(\int_{B_{r_0}}  V_2 (u^*)^{q} dx\Big)^{\frac2{q}}\geq  c c_q^{2/q}\Big(\int_{B_{r_0}}  V_2 u^{q} dx\Big)^{\frac2{q}}.
      \end{eqnarray*}  
The proof is complete. 
\hfill$\Box$\medskip 
 
\section{Embedding results for $\hat\cH_0^1(B_1)$ and related spaces}
\subsection{Embedding and compactness} 

We prove Theorems \ref{embedding} and \ref{embedding-t} here.

\begin{proof}[{\bf Proof of Theorem \ref{embedding}}]
 $(i)$ 
  Recall that the Sobolev   space $W^{1,q}_0(B_1)$ is  the completion of $C_c^{\infty}(B_1)$ under the norm
$$\|u\|_{1,q}=\Big(\int_{B_1}  |\nabla u |^q dx\Big)^{\frac1q}. $$
We first show that 
$$  \hat\cH_0^1(B_1)\hookrightarrow  W^{1,q}_0(B_1)\quad{\rm for\ any\ } q\in[1,2). $$

For any given radially symmetry function $u_0\in C_c^{\infty}(B_1)$, we set 
$$w_0(x)=(-\ln |x|)^{-\frac12} u_0(x).$$
By direct computations,
\begin{eqnarray*} 
   \int_{B_1} \Big(|\nabla u_0|^2-\frac14\frac{u_0^2}{|x|^2(-\ln |x|)^2}\Big)   dx
    =   \int_{B_1}  |\nabla w_0|^2  (-\ln |x|) dx 
    \end{eqnarray*}
and 
 \begin{eqnarray*} 
 \int_{B_1} |\nabla u_0|^qdx &=&  \int_{B_1}  \Big|(-\ln |x|)^{\frac12}  \nabla w_0-\frac12(-\ln|x|)^{-\frac12} \frac{x}{|x|^2}w_0 \Big|^qdx
 \\[2mm]&\leq &2^q \int_{B_1} (-\ln |x|)^{\frac q2} |\nabla w_0|^qdx+\int_{B_1} (-\ln|x|)^{-\frac q2}|x|^{-q} |w_0|^qdx.
 \end{eqnarray*} 
Since
 \begin{eqnarray*}  
 \int_{B_1} (-\ln |x|)^{\frac q2} |\nabla w_0|^qdx\leq \Big(\int_{B_1}    (-\ln |x|) |\nabla w_0|^2 dx\Big)^{\frac q2} |B_1|^{1-\frac q2},
 \end{eqnarray*} 
  by the H\"older inequality and (\ref{eqn change1})
  \begin{eqnarray*}  
 &&\int_{B_1}  (-\ln|x|)^{-\frac q2}|x|^{-q} |w_0|^q dx\\[2mm]&\leq& \Big(\int_{B_1}  |x|^{-2} (-\ln|x|)^{-1}\big(1+|\ln \ln \frac{1}{|x|}|\big)^{-2}  w_0^2  dx\Big)^{\frac q2}  \Big(\int_{B_1} \big(1+|\ln \ln \frac{1}{|x|}|\big)^{\frac{2q}{2-q}} dx\Big)^{1-\frac q2}
\\[2mm]&\leq& c\Big(\int_{B_1}   (-\ln |x|) |\nabla w_0|^2 dx\Big)^{\frac q2} ,
 \end{eqnarray*} 
 where we have used $\int_{B_1} \big(1+|\ln \ln \frac{1}{|x|}|\big)^{\frac{2q}{2-q}} dx<+\infty$.

As a consequence, we obtain that for $q\in[1,2)$, there exists $c=c(q)$ such that 
$$ \int_{B_1} |\nabla u_0|^qdx\leq  c\Big(\int_{B_1}   (-\ln |x|) |\nabla w_0|^2 dx\Big)^{\frac q2}= c\Big(\int_{B_1} \Big(|\nabla u_0|^2-\frac14\frac{u_0^2}{|x|^2(-\ln |x|)^2}\Big)   dx\Big)^{\frac q2}.$$
Hence the embedding 
  $\hat\cH_0^1(B_1)\hookrightarrow W^{1,q}_0(B_1)$ is continuous for $q\in[1,2)$.\vs
  
   $(ii)$ Since 
  $W^{1,q}_0(B_1)\hookrightarrow L^p(B_1)$ is compact for $p<\frac{2q}{2-q}$ (see \cite[Theorem 7.10]{GT}),
  and $\frac{2q}{2-q}\to+\infty$ as $q\to 2^-$, for any $p\in [1,+\infty)$, we may take $q<2$ close to 2 such that $p< \frac{2q}{2-q}$, and then the embedding
  $\hat\cH_0^1(B_1)\hookrightarrow W^{1,q}_0(B_1)$ is continuous, and the embedding $W^{1,q}_0(B_1)\hookrightarrow L^p(B_1)$ is compact. It follows that
  $\hat\cH_0^1(B_1)\hookrightarrow L^p(B_1)$ is compact. This proves the conclusion in $(ii)$ with $\alpha=0$.

      When $\alpha\in(0,2)$,   let 
$$\theta=\frac{2+\alpha}{2\alpha}>1,\qquad \theta'=\frac{\theta}{\theta-1}=\frac{2+\alpha}{2-\alpha}.$$ 
Then
by the H\"older inequality
  \begin{eqnarray*}  
\big( \int_{B_1} u^p |x|^{-\alpha} dx\big)^{\frac1p} &\leq& \big(\int_{B_1}  u^{p\theta'}   dx\big)^{\frac 1{p\theta'}} \big( \int_{B_1} |x|^{-\alpha \theta} dx\big)^{\frac1\theta},
 \end{eqnarray*}
 where $\alpha \theta\in (1,2)$.

 Now the compactness of embedding follows from the fact that 
  $W^{1,q}_0(B_1)\hookrightarrow L^{p\theta'}(B_1)$ is compact for 
  $$p\theta'<\frac{2q}{2-q},$$
  which is satisfied for any given $p\geq 1$ if we choose 
    $q\in(0,2)$ close enough to 2. 
    \end{proof}

    \noindent{\bf Proof of Theorem \ref{embedding-t}. }    Let
 $u_n \rightharpoonup 0$ weakly in $\hat\cH^1_0(B_1)$ as $n\to+\infty$. Then $\{ u_n\}$ is uniformly bounded in $\hat\cH^1_0(B_1)$. For given $\epsilon>0$, there exists $\sigma_0\in(0,\frac18)$ such that for $x\in B_{\sigma_0}\cup \big(B_1\setminus B_{1-\sigma_0}\big)$ there holds
$$ V(x) |x|^{2} (-\ln |x|)^{2} \big(1+|\ln \ln \frac{1}{|x|}|\big)^{2}  \leq  2V(x) |x|^{2} (-\ln |x|)^{2} \big|\ln \ln \frac{1}{|x|}\big|^{2}\leq \epsilon.  $$ 
 Thus, by the improved Leray inequality (\ref{e 3.1-2}) we obtain 
\begin{eqnarray*} 
 && \int_{ B_{\sigma_0}\cup \big(B_1\setminus B_{1-\sigma_0}\big)} |u_n(x)|^2V dx\\
  &\leq&  \epsilon  \int_{ B_{\sigma_0}\cup \big(B_1\setminus B_{1-\sigma_0}\big)} |u_n(x)|^2 |x|^{-2} (-\ln |x|)^{-2} \big(1+|\ln \ln \frac{1}{|x|}|\big)^{-2}  dx 
 \\[2mm]&\leq&  \epsilon  \int_{ B_1} |u_n(x)|^2 |x|^{-2} (-\ln |x|)^{-2} \big(1+|\ln \ln \frac{1}{|x|}|\big)^{-2}  dx 
 \\[2mm]& \leq&  \epsilon  \|u_n\|_{\hat\cH^1_0}^2 .
 \end{eqnarray*}  
Moreover, by the compactness embedding $\hat\cH^1_0(B_1) \hookrightarrow L^2(B_1)$,
 \begin{eqnarray*} 
  \int_{B_{1-\sigma_0}\setminus B_{\sigma_0}} |u_n(x)|^2V dx &\leq&  \Big( \max_{\overline{B_{1-\sigma_0}\setminus B_{\sigma_0} }}V \Big) \int_{B_1} |u_n(x)|^2 dx 
  \leq \epsilon \mbox{  for all large $n$. }
 \end{eqnarray*}

Therefore  we obtain
$$
	 \int_{B_1} |u_n(x)|^2  V dx \to 0\quad{\rm as}\ n\to+\infty.
$$
It follows that  $\hat\cH^1_0(B_1)$ is compactly embedded in $L^2(B_1, V dx)$. 

To obtain the compactness of the embedding 
$$\hat\cH^1_0(B_{r_q})\hookrightarrow  L^q\big(B_{r_q},  Vdx\big),$$
we use \eqref{e 3.1-q ral} instead of (\ref{e 3.1-2}), and repeat the above analysis with obvious modifications.
The proof is completed. \hfill$\Box$

\subsection{General weighted spaces  } 

Recall that  $\cH^1_{\mu,0}(B_1)$ is  the completion of $C_c^\infty(B_1)$ under the norm
$$\|u\|_\mu=\sqrt{\int_{B_1}\Big(   |\nabla u |^2  +\mu \frac{ u^2}{|x|^2(-\ln|x|)^2} \Big)dx} $$
and it is a Hilbert space with the inner product 
$$\langle u,v\rangle_\mu=\int_{B_1}  \Big( \nabla u \cdot\nabla v\, +\mu\frac{ uv}{|x|^2(-\ln|x|)^2}\Big) dx.$$ 
Recall also  that  
$$
\cH_0^1(B_1)=\cH^1_{0, 0}(B_1) \quad{\rm and}\quad \hat\cH^1_0(B_1)=\cH^1_{-\frac14, 0}(B_1).$$ 

We have the following result.
 
 \begin{theorem}\label{th:main-1} The following conclusions hold:
 
$(i)$ For $\mu>-\frac14$,  $\cH^1_{\mu, 0}(B_1)= \cH^1_{0}(B_1)$. \smallskip
 
$(ii)$  $ \cH^1_{0}(B_1)\varsubsetneqq \hat\cH_{0}^1(B_1)$.	  
\end{theorem}
 \noindent{\bf Proof.}  $(i)$ By (\ref{ei 1.1}),  for $u\in \cH^1_{\mu,0}(B_1)$, we have
\begin{eqnarray*} 
  \|u\|_\mu^2  =   \int_{B_1}   |\nabla u |^2dx  +\mu \int_{B_1}\frac{u^2}{|x|^2(-\ln|x|)^2} dx
 \geq   (1+4 \min\{0,\mu\} ) \int_{B_1}   |\nabla u |^2dx,
    \end{eqnarray*} 
\begin{eqnarray*} 
  \|u\|_\mu^2   \leq   \big(1+4\max\{0,\mu\}\big) \int_{B_1}   |\nabla u |^2dx.
    \end{eqnarray*} 
The desired conclusion then follow immediately.

$(ii)$ If $u\in \cH^1_{0}(B_1)$, then by
 Leray's inequality, 
$$\int_{B_1}\frac{u^2}{|x|^2(-\ln|x|)^2} dx<+\infty.$$ Hence $u\in \hat\cH^1_0(B_1)$, which proves
 $ \cH^1_{0}(B_1)\subset  \hat\cH^1_{ 0}(B_1)$.
On the other hand, let  $h_0(x)=(-\ln |x|)^{\frac12} \eta_0(2|x|)$, where
 $\eta_0:\R\to[0,1]$ is a smooth cutoff function such that 
$$
  \eta_0(t) =\left\{ \arraycolsep=1pt
\begin{array}{lll}
 1\ \ \ \ &{\rm for}\ \,   |t|<\frac12,
\\[2mm] 
0\ \  &{\rm for}\ \,   |t|>1 . 
\end{array}
\right.  
$$ 
Direct computation shows that $h_0\in \hat\cH^1_{0}(B_1)$, but it is not in $\cH^1_{0}(B_1)$.
\hfill$\Box$\medskip

Recall that $\cH^1_{_V,\mu, 0}(B_1)$ with $\mu\geq-\frac14$ is  the completion of $C_c^\infty(B_1\setminus\{0\})$ under the norm 
 $\|\cdot\|_{_V,\mu}$ defined in (\ref{norm V}).

 \begin{corollary}\label{cr:main-0}
Let $\Omega$ be a  domain in $B_1$ containing the origin, $V:B_1\setminus\{0\}\to [0,+\infty)$ be  a continuous function satisfying
 \begin{equation}\label{con V-l0}\begin{cases}
\displaystyle \limsup_{|x|\to0^+} V(x) |x|^2(-\ln |x|)^2<+\infty,\\[2.5mm]
\displaystyle \limsup_{|x|\to1^-} V(x)(1-|x|)^2<+\infty.
\end{cases}
\end{equation} 
Then for any $\mu>0$,     $\cH^1_{_V,\mu, 0}(\Omega)= \cH^1_{0}(\Omega)$.  
 \end{corollary}
\noindent{\bf Proof. } It follows from (\ref{con V-l0})  
that   
$$   V(x)\leq c |x|^{-2}(-\ln |x|)^{-2} \ \mbox{ for some $c>0$ and all $x\in B_1\setminus\{0\}$},$$
which implies, in view of $\mu>0$,
$$\cH^1_{_V,\mu, 0}(\Omega)\subset \cH^1_{c\mu,0}(\Omega)=\cH^1_{0}(\Omega).$$
On the other hand,    Leray's inequality implies $\cH^1_0(\Omega)\subset \cH^1_{_V,\mu, 0}(\Omega)$.\hfill$\Box$\medskip

 \begin{corollary}\label{cr:main-1}
Let $\Omega$ be a  domain in $B_1$ containing the origin, $V:B_1\setminus\{0\}\to [0,+\infty) $  be  continuous  and 
 verify
 \begin{equation}\label{con V-nonr}
 V(x)\leq \frac{1}{|x|^2(-\ln |x|)^2}.
\end{equation} 
 Then the following conclusions hold:
 
$(i)$ For $\mu\in(-\frac14,0)$,  $\cH^1_{_V,\mu, 0}(\Omega)= \cH^1_{0}(\Omega)$. \smallskip
 
$(ii)$ For $\mu=-\frac14$, assume additionally that \eqref{con V-l2} holds,
 then $ \cH^1_{0}(\Omega)\varsubsetneqq \cH^1_{_V,-1/4, 0}(\Omega)$.	  
\end{corollary}
{\bf Proof. }    For $\mu\in(-\frac14, 0)$, by  (\ref{con V-nonr}) and the Leray inequality, we easily see that $\cH^1_{_V,\mu, 0}(\Omega)= \cH^1_{0}(\Omega)$ and so $(i)$ holds true. 

{$(ii)$}  Let $r_0\in(0,1)$ be such that $B_{r_0}(0)\subset \Omega$,  and then take  $h_0(x):=(-\ln |x|)^{\frac12} \eta_0(\frac {2}{r_0}|x|)$; then  $h_0\not\in\cH^1_{0}(B_{r_0})$ by a direct computation. However $h_0\in \cH^1_{_V,-1/4, 0}(B_{r_0})$ since $h_0\in\hat\cH_0^1(B_{r_0})$ and
$$ \int_{B_{r_0}} \Big|V-\frac1{|x|^2(-\ln|x|)^2} \Big| \,h_0^2\, dx<+\infty,$$
which is guaranteed by (\ref{con V-l2}).\hfill$\Box$ \medskip

\medskip

\noindent{\bf Proof of Theorem \ref{coro 2.1}. } 
$(i)$ Let  $w_0\in C_c^{1}(B_1)$ and  
$$u_0(x)=(-\ln |x|)^{\frac12} w_0(x). $$
Then $u_0\in C^1(B_1\setminus \{0\})\cap   \hat\cH_{0}^1(B_1)$.

Direct computation shows  that 
\begin{eqnarray*} 
  \cI(u_0) &=&  \int_{B_1} \Big(|\nabla u_0|^2-\frac14\frac{u_0^2}{|x|^2(-\ln |x|)^2}\Big)   dx
   \\[2mm]&=& \int_{B_1} \Big(\Big|\frac12(-\ln|x|)^{-\frac12} \frac{x}{|x|^2}w_0 +(-\ln |x|)^{\frac12}  \nabla w_0\Big|^2-\frac14\frac{w_0^2}{|x|^2(-\ln |x|)}\Big)   dx
   \\[2mm]&=&  \int_{B_1}  |\nabla w_0|^2  (-\ln |x|) dx+\frac12\int_{B_1}\frac{x\cdot \nabla(w_0^2)}{|x|^2}dx
     \\[2mm]&=&  \int_{B_1}  |\nabla w_0|^2  (-\ln |x|) dx,
    \end{eqnarray*}
    where 
 $$\int_{B_1}\frac{x\cdot \nabla(w_0^2)}{|x|^2}dx=\int_{B_1} {\rm div}(\frac{xw_0^2}{|x|^2})dx=0. $$
Then (\ref{eqn change1}) and (\ref{eqn change1-q}) follow from the embedding inequality (\ref{e 3.1-2}) in Remark \ref{rm em1} by  setting
 $w_0(x)=(-\ln |x|)^{-\frac12} u_0$. \vs

Part $(ii)$  (a) and (b) follow from the proofs of Theorem \ref{embedding} and Theorem \ref{embedding-t} respectively. 
\hfill$\Box$ \vs

     \setcounter{equation}{0}
     \section{ Trudinger-Moser type inequalities for radial functions }
  
  Define
$$\cH^1_{{\rm rad},\mu, 0}(B_1):=\left\{w\in \cH^1_{\mu, 0}(B_1):\, w\ {\rm is\ radially\ symmetric}\right\}.$$
We will prove the following two theorems in this section.

 \begin{theorem}\label{teo 1} Let $\mu>-\frac14$. Then
 \[
\sup_{u\in\mathcal H^1_{{\rm rad}, \mu, 0}(B_1),\ \|u\|_{\mu}\leq 1}\int_{B_1}e^{4 \pi     \sqrt{1+4\mu} \; u^2 }dx<\infty,
 \]
 and the result fails when $4 \pi     \sqrt{1+4\mu}$ is replaced by any $\alpha>4 \pi     \sqrt{1+4\mu}$. 
  \end{theorem}

Note that in the above result $4 \pi     \sqrt{1+4\mu}=m_\mu\to0$ as $\mu\to -\frac14$,  
which  suggests that the  inequality should be different in the case $\mu=-\frac14$. Our Trudinger-Moser type inequality for   $\mu=-\frac14$ is the following. 
 
  \begin{theorem}\label{teo 2}
   $(i)$  For any $p\in(0,1)$ and any  $\alpha>0$, there exists $c=c_{p,\alpha}$ depending on $p$ and $\alpha$ such that  for every $u\in \cH^1_{{\rm rad}, {-1/4}, 0}(B_1)$ with $\|u\|_{-1/4}\leq 1$,
 there holds
$$\int_{B_1} e^{\alpha |u|^p}dx\leq  c_{p,\alpha}.$$
  
  $(ii)$ For any $p\geq1$ and any $\alpha>0$,
   there exists a sequence $\{u_n\}$ such that  $\|u_n\|_{-1/4}\leq 1$,
$$\int_{B_1} e^{\alpha |u_n|^p }dx\to +\infty\quad{\rm as}\ \, n\to+\infty.$$
 \end{theorem}
 \medskip
 
 The following simple lemma will be very useful later.

 \begin{lemma}\label{lm 2.1-p-1}
Let $\mu>-\frac14$,   $t_1>2$,  $\tau_0=\frac{1-\sqrt{1+4\mu}}{2}$,
and  
$$\cJ_{t_1}(w):= \frac1{ 1-2\tau_0 } \int_2^{+\infty} [w'(t)]^2 t^{2\tau_0} dt \ \mbox{ for } \ w\in \bX_{t_1}$$
with 
$$ \begin{cases}\bX_{t_1}:=\Big\{ v\in C^{0,1} ([2,\infty)): v(2)=0,\ v(t_1)=C_{t_1}( t_1^{1-2\tau_0} -2^{1-2\tau_0})\Big\},\\[2mm]
 C_{t_1}:=\big(t_1^{1-2\tau_0}-2^{1-2\tau_0}\big)^{-1/2}. 
 \end{cases}$$
 Then   $\cJ_{t_1}(w)\geq 1$, with equality holding only if $w=\zeta_{t_1}$, where
\begin{eqnarray}  \label{cla f1}
\zeta_{t_1}(t):= 
 \left\{ \arraycolsep=1pt
\begin{array}{lll}
  C_{t_1}(t^{1-2\tau_0}-2^{1-2\tau})\ \ \ \ &{\rm for}\ \,   t\in [2,t_1],
\\[2mm]
 C_{t_1}(t_1^{1-2\tau_0} -2^{1-2\tau}) \ \  &{\rm for}\ \,   t\in(t_1,+\infty).
\end{array}
\right.  
\end{eqnarray} 

 \end{lemma}
\noindent{\bf Proof.} Direct computation shows that 
\begin{eqnarray*} 
\cJ_{t_1}(\zeta_{t_1})&=&  \frac1{ 1-2\tau_0 }  \int_2^{t_1} [\zeta_{t_1}'(t)]^2 t^{2\tau_0} dt 
\\&=&  \big(t_1^{1-2\tau_0}-2^{1-2\tau_0}\big)^{-1} (1-2\tau_0 )\int_2^{t_1} t^{-2\tau_0}dt   \,=1.  
\end{eqnarray*} 
Given $w\in \bX_{t_1}$,  let $v:=w-\zeta_{t_1}$; then 
$v(2)=v(t_1)=0$ and
we have  
\begin{eqnarray*} 
\cJ_{t_1}(w)&=&  \frac1{ 1-2\tau_0 }  \int_2^{+\infty} \big(v'+\zeta_{t_1}'\big)^2 t^{2\tau_0} dt
\\&=&    \frac1{ 1-2\tau_0 }   \int_2^{+\infty}  \big[v'(t)\big]^2 t^{2\tau_0} dt +2 C_{ t_1}  \int_2^{t_1}  v'(t) dt  +1 
\\&=&    \frac1{ 1-2\tau_0 }  \int_2^{+\infty} \big[ v'(t)\big]^2 t^{2\tau_0} dt +1 
\\&\geq &1, 
\end{eqnarray*} 
where the last inequality becomes an equality only if $\int_2^{+\infty}  v'(t)^2 t^{2\tau_0} dt=0$, which implies $v'(t)=0$ in $(2,+\infty)$, and we thus obtain, in view of $v(2)=0$,  $v(t)=0$  in $(2,+\infty)$. 
  \hfill$\Box$\medskip
 
   \subsection{Proof of Theorem \ref{teo 1}}
  Let 
 $$C^\infty_{{\rm rad},c}(B_1):=\big\{w\in C^\infty_c(B_1)\!: \, w\ {\rm is\ radially\ symmetric}  \big \}.$$
Note that $C^\infty_{{\rm rad},c}(B_1)$ is dense in $\cH^1_{{\rm rad},\mu, 0}(B_1)$.

 \begin{lemma}\label{teo 3.1}
Suppose $\mu>-\frac14$ and  $m_\mu =4\pi  \sqrt{1+4\mu}$. 
Then    
$$
\sup_{u\in \cH^1_{{\rm rad},\mu, 0}(B_1),\, \|u\|_\mu\leq 1}\int_{B_1} e^{m_\mu u^2}dx<\infty.$$

 \end{lemma}
 \noindent {\bf Proof. } Since   $C^\infty_{{\rm rad},c}(B_1\setminus\{0\})$  is dense in $\cH^1_{{\rm rad},\mu, 0}(B_1)$,   to prove 
  the desired inequality,  we start with estimating $\int_{B_1} e^{m_\mu u_0^2}dx$ for   an arbitrary  function $u_0\in C^\infty_{{\rm rad},c}(B_1\setminus\{0\})$
  with $\|u_0\|_\mu\leq 1$.   
 
 Denote 
 \[
 \tau_0=\tau_0(\mu):=\frac{1-\sqrt{1+4\mu}}{2},
 \]
 and write   
$$u_0(r)=(1-\ln r)^{\tau_0} v_0(r)\quad{\rm for}\ \, r\in(0,1).$$
Then $v_0(r)$ has compact support in $(0,1)$.

Direct computations give 
\begin{eqnarray*} 
  \cI(u_0) &=& 2\pi \int^1_0 \Big(u_0'(r)^2+\mu \frac{u_0(r)^2}{r^2(1-\ln r)^2}\Big)  rdr
 \\[2mm] &=&2\pi \int^1_0  v_0'(r)^2  (1-\ln r)^{2\tau_0} rdr - 4\pi  \tau_0  \int^1_0 v_0'(r)v_0(r)(1-\ln r)^{2\tau_0-1}dr 
 \\&&+2\pi(\tau_0^2+\mu)\int^1_0  v_0^2(r)(1-\ln r)^{2\tau_0-2} r^{-1}dr
 \\[2mm] &=&2\pi \int^1_0  v_0'(r)^2  (1-\ln r)^{2\tau_0} rdr  
 -2\pi(\tau_0^2-\tau_0-\mu)\int^1_0  v_0^2(r)(1-\ln r)^{2\tau_0-2} r^{-1}dr
 \\[2mm]&=&  2\pi \int^1_0  v_0'(r)^2  (1-\ln r)^{2\tau_0} rdr,   
    \end{eqnarray*} 
 where we have used $\tau_0^2-\tau_0-\mu=0$. Note also that $\tau_0\in(0,\frac12)$ for $\mu\in(-\frac14,0)$ and $\tau_0<0$ for $\mu>0$.

  Let 
  $$  r=e^{1-\frac12 t},\qquad v_0(r)=w_0(t).$$
  Then $w_0(t)=0$ for $t>2$ close to 2, and $w(\infty)=0$. Moreover,
\[
 \cI(u_0)=  2\pi \int^1_0  v_0'(r)^2  (1-\ln r)^{2\tau_0} rdr=4^{1-\tau_0}\pi\int_2^\infty w_0'(t)^2 t^{2\tau_0}dt.
  \]
Thus, for given $t>2$, 
    \begin{eqnarray*} 
  w_0^2(t)    =   \Big( \int^t_2  w_0'(s)  ds\Big)^2 &\leq&  \int^t_2  w_0'(s)^2   s^{ 2 \tau_0 }  ds    \int^t_2     s^{ -2 \tau_0 } ds  
  \\[2mm]&=  & (1-2\tau_0)^{-1} \int^t_2  w_0'(s)^2   s^{ 2 \tau_0 }  ds \, ( t^{1- 2\tau_0 }-2^{1-2\tau_0})
    \end{eqnarray*}
 and       
 \begin{eqnarray*} 
  1\geq \cI(u_0)    =   4^{1-\tau_0}\pi  \int^\infty_2  w_0'(s)^2 s^{2\tau_0 }  ds
 \geq 4^{1-\tau_0}\pi (1-2\tau_0) w_0^2(t)  (t^{1-2\tau_0}-2^{1-2\tau_0})^{-1} ,
    \end{eqnarray*}
which implies that 
   \begin{eqnarray} 
  w_0^2(t)   \leq    \frac{  t^{1- 2\tau_0 }-2^{1-2\tau_0}}{ 4^{1-\tau_0}\pi (1-2\tau_0)} \quad \mbox{ for } t>2. \label{eqvi 1-2.1}
    \end{eqnarray}

  We now have, by simple calculations,
  \begin{eqnarray*} 
 \cN(u_0) :=\int_{B_1} e^{m_\mu u_0^2(|x|)}dx&=&2\pi \int^1_0  e^{m_\mu u_0^2(r)}r dr\\
  &=& \pi e^2  \int^{+\infty}_2  e^{ m_\mu 4^{-\tau_0} s^{2\tau_0}w^2_0(s)-s }  ds \\
  &=&\pi  e^2 \int^{+\infty}_2  e^{   s^{2\tau_0} \tilde w^2_0(s)-s }  ds,
     \end{eqnarray*} 
  where
  \[\tilde w_0(t):=2^{-\tau_0}\sqrt{m_\mu}\, w_0(t).
  \]
  We note that
  \[
  4^{1-\tau_0}\pi (1-2\tau_0)=4^{-\tau_0}m_\mu 
  \]
  and hence, by
     (\ref{eqvi 1-2.1}), 
     \[
     \tilde w_0^2(t)\leq t^{1- 2\tau_0 }-2^{1-2\tau_0},\quad t^{2\tau_0}\tilde w_0^2(t)-t\leq -2^{1-2\tau_0}t^{2\tau_0}.
     \]
    Therefore, when $\mu\in (-1/4, 0)$ and hence $\tau_0\in (0, 1/2)$,
     \[
      \cN(u_0)=\pi  e^2 \int^{+\infty}_2  e^{   s^{2\tau_0} \tilde w^2_0(s)-s }  ds\leq \pi e^2\int_2^\infty e^{-2^{1-2\tau_0}s^{2\tau_0}}ds<\infty.
     \] 
This completes the proof for the case $\mu\in (-1/4, 0)$. \medskip

It remains to consider the case $\mu>0$ and hence $\tau_0<0$.
Since
     \[
     \tilde w_0^2(t)\leq t^{1- 2\tau_0 }-2^{1-2\tau_0}=\big[C_t\, (t^{1-2\tau_0}-2^{1-2\tau_0})\big]^2 \ \mbox{ for } t>2,
     \]
     we have
  $$
  \sigma_0:=1-\sup_{t>2}\frac{\tilde w_0(t)}{ C_t\, (t^{1- 2\tau_0 }-2^{1-2\tau_0})} \geq0. 
  $$
As $w_0(t)$ and hence $\tilde w_0(t)=0$  for $t>2$ close to 2, and
\[\mbox{ $\displaystyle w_0(\infty)=\tilde w_0(\infty)=0$, $\displaystyle\lim_{t\to\infty}C_t (t^{1- 2\tau_0 }-2^{1-2\tau_0})=\infty$,}
\]
there exists $t_1\in(2,+\infty)$ such that 
 $$\sigma_0=1- \frac{\tilde w_0(t_1)}{C_{t_1} (t_1^{1- 2\tau_0}-2^{1-2\tau_0})}.$$
 Let us also note that
   \begin{eqnarray*} 
 \cJ_{t_1}(\tilde w_0)&=& \frac1{ 1-2\tau_0 } \int_2^{+\infty} [\tilde w_0'(t)]^2 t^{2\tau_0} dt
 \\[1mm]&  =&4^{1-\tau_0}\pi  \int^\infty_2  [w_0'(t)]^2 t^{2\tau_0 }\,  dt\ \,=  \cI(u_0)\leq 1.
 \end{eqnarray*} 
 
 We first show that $\sigma_0=0$ cannot happen. Indeed, 
 if $\sigma_0=0$,  then
 $$\tilde w_0(t_1) =C_{t_1}(t_1^{1- 2\tau_0}-2^{1-2\tau_0}),$$ 
 and  it follows from Lemma \ref{lm 2.1-p-1}  that 
$$\tilde w_0\equiv \zeta_{t_1}\quad {\rm in}\ \, (2,\infty), $$ 
  where  $\zeta_{t_1}$ is defined in (\ref{cla f1}). Thus $\tilde w_0(t)>0$ for all $t>2$, which contradicts the fact that $\tilde w_0(t)=0$ for $t>2$ close to 2.
  Thus we must have $\sigma_0\in (0,1)$.

 We now have
 \[
 s^{2\tau_0}\tilde w^2_0(s)-s\leq (1-\sigma_0)^2C^2_ss^{2-2\tau_0}-s=\big[(1-\sigma_0)^2C_s^2 s^{1-2\tau_0}-1\big]s \ \ \mbox{ for }\, s>2,
 \]
 and also, by
     (\ref{eqvi 1-2.1}), 
     \[
      s^{2\tau_0}\tilde w_0^2(s)-s\leq -(1-2^{2\tau_0-2}\pi^{-1})s \ 
      \mbox{ for } s>2.
 \]
 Therefore, 
   \begin{eqnarray*} 
 \cN(u_0)  
 &=& \pi e^2 \int^ {3}_2  e^{s^{2\tau_0}\tilde w^2_0(s)-s } ds +  \pi e^2 \int^\infty_{3 }  e^{s^{2\tau_0}\tilde w^2_0(s)-s}  ds\\
 &\leq& \pi e^2 \int^ {3}_2  e^{ -(1-2^{2\tau_0-2}\pi^{-1}) s } ds+\pi e^2 \int^\infty_{3 }  e^{s^{2\tau_0}\tilde w^2_0(s)-s}  ds.
       \end{eqnarray*} 
  Define     
       \[
 G(\sigma_0):= \int^\infty_{3 }  e^{s^{2\tau_0}\tilde w^2_0(s)-s}  ds.
 \]
 Clearly,
 \[
 (1-\sigma_0)^2C_s^2 s^{1-2\tau_0}-1=(1-\sigma_0)^2\frac{s^{1-2\tau_0}}{s^{1-2\tau_0}-2^{1-2\tau_0}}-1<-\sigma_0 \ \mbox{ for all large } s.
 \]
 It follows that 
 \[
 G(\sigma_0)\leq \tilde G(\sigma_0):=\int_ {3}^\infty  e^{ [(1-\sigma_0)^2C_s^2 s^{1-2\tau_0}-1]s} ds<\infty.
 \]
 Since $\tilde G(\sigma_0)$ is a continuous decreasing function of $\sigma_0$ for $\sigma_0\in (0,1)$, it is now clear that
 to obtain our desired result, it suffices to show that
 \[
 \limsup_{\sigma_0\to 0^+}G(\sigma_0)<\infty.
 \]
 
 Thanks to 
\[ \begin{aligned}
\cJ_{1}(\zeta_{t_1})&=  \frac1{ 1-2\tau_0 }  \int_2^{t_1} \zeta_{ t_1}'(s)^2 s^{2\tau_0} ds =1,  \\
 \cJ_{t_1}(\tilde w_0)&= \frac1{ 1-2\tau_0 } \int_2^{\infty} [\tilde w_0'(t)]^2 t^{2\tau_0} dt=\cI(u_0)\leq 1
 \end{aligned}
\]
and 
\[
\int_2^{t_1} \big[\tilde w_0'(s)- (1-\sigma_0)\zeta_{t_1}'(s)\big] \zeta_{t_1}'(s) s^{2\tau_0}ds=C_{t_1}(1-2\tau_0)\int_2^{t_1} \big[\tilde w_0'(s)-(1-\sigma_0) \zeta_{t_1}'(s)\big]ds=0,
\]
  we obtain 
   \begin{align*}
   &\int_2^{t_1} \big(\tilde w_0'-(1-\sigma_0) \zeta_{t_1}'\big)^2 s^{2\tau_0}ds +
   \int_{t_1}^\infty (\tilde w_0')^2 s^{2\tau_0} ds\\[1mm]
   =&\ \int_2^{t_1} \big[\tilde w_0'-(1-\sigma_0) \zeta_{t_1}'+(1-\sigma_0)\zeta_{t_1}' \big]^2 s^{2\tau_0}ds-\int_2^{t_1} \big[(1-\sigma_0)\zeta_{t_1}' \big]^2 s^{2\tau_0}ds +
   \int_{t_1}^\infty (\tilde w_0')^2 s^{2\tau_0} ds   \\[1mm]
   \leq &\ \, (1-2\tau_0)[1-(1-\sigma_0)^2] \\[1mm]
   = &\ \,  (1-2\tau_0)\sigma_0(2-\sigma_0). 
   \end{align*}
 Therefore, for $s>t_1$, 
 \begin{eqnarray*}  
 \tilde w_0(s) &\leq& \tilde w_0(t_1)+ \left(\int_{t_1}^s\tilde w_0'(t)^2 t^{2\tau_0} dt\right)^{\frac12}\Big(\int_{t_1}^s  t^{-2\tau_0} dt\Big)^{\frac12}
 \\[2mm]&\leq& \min\Big\{ t_1^{\frac{1-2\tau_0}2}, (1-\sigma_0)C_{t_1}t_1^{1-2\tau_0}\Big\}+\sqrt{\sigma_0(2-\sigma_0) }\Big(s^{\frac{1-2\tau_0}2}-t_1^{\frac{1-2\tau_0}2}\Big),
  \end{eqnarray*} 
  and for $s\in(2,t_1]$,
 \begin{eqnarray*}  
 \tilde w_0(s)-(1-\sigma_0) \zeta_{t_1}(s)&\leq&   \left(\int_2^s\big[\tilde w_0'-(1-\sigma_0) \zeta_{t_1}'\big]^2 t^{2\tau_0} dt\right)^{\frac12}\Big(\int_2^s  t^{-2\tau_0} dt\Big)^{\frac12}
 \\[2mm]&\leq&  \sqrt{\sigma_0 (2-\sigma_0)} s^{\frac{1-2\tau_0}2},
  \end{eqnarray*} 
\begin{eqnarray*}  
 \tilde w_0(s)-(1-\sigma_0) \zeta_{t_1}(s)&\leq&   \left(\int_s^{t_1}\big[\tilde w_0'-(1-\sigma_0) \zeta_{t_1}'\big]^2 t^{2\tau_0} dt\right)^{\frac12}\Big(\int_s^{t_1}  t^{-2\tau_0} dt\Big)^{\frac12}
 \\[2mm]&\leq&  \sqrt{\sigma_0 (2-\sigma_0)}\Big(t_1^{\frac{1-2\tau_0}2}- s^{\frac{1-2\tau_0}2} \Big).
  \end{eqnarray*} 
  Thus
  \[
  s^{2\tau_0}\tilde w_0(s)^2\!-\!s\!\leq\!\begin{cases}\big[\min\Big\{ t_1^{\frac{1\!-\!2\tau_0}2}\!\!, (1\!-\!\sigma_0)C_{t_1}t_1^{1\!-\!2\tau_0}\Big\}s^{\tau_0}\!+\!\sqrt{\sigma_0(2\!-\!\sigma_0)} (s^{\frac 12}\!-\!s^{\tau_0}t_1^{\frac 12 -\tau_0})\big]^2\!-\!s, &\!\! s\!>\!t_1,\\[2mm]
  \big[(1\!-\!\sigma_0)s^{\tau_0}\zeta_{t_1}(s)\!+\!\sqrt{\sigma_0(2-\sigma_0)}s^{\tau_0}\min\big\{ s^{\frac{1-2\tau_0}2}\!, t_1^{\frac{1\!-\!2\tau_0}2}\!-\! s^{\frac{1-2\tau_0}2} \big\}\big]^2\!-\!s,&\!\! s\!\in\! (2, t_1].
  \end{cases}
  \]
  
  If $t_1\leq 3$, then 
  \[
  s^{2\tau_0}\tilde w_0(s)^2-s\leq\big[\sqrt 3+\sqrt{\sigma_0(2-\sigma_0) s}\big]^2-s \mbox{ for } s>3,
  \]
  and hence
  \[
  \limsup_{\sigma_0\to 0^+}G(\sigma_0)\leq \lim_{\sigma_0\to 0^+} \int_3^\infty e^{\big[\sqrt 3+\sqrt{\sigma_0(2-\sigma_0) s}\big]^2-s}dx<\infty.
  \]
  
  If $t_1\in (3, M]$ for some $M>0$, then
  \begin{align*}
  G(\sigma_0)&\leq \int^{t_1}_3 e^{\big[(1-\sigma_0)s^{\tau_0}\zeta_{t_1}(s)+\sqrt{\sigma_0(2-\sigma_0)s}\big]^2-s}ds 
  +\int^\infty_{t_1}e^{\big[(1-\sigma_0)C_{t_1}t_1^{1-2\tau_0}s^{\tau_0}+\sqrt{\sigma_0(2-\sigma_0) s}\big]^2-s}ds\\
  &\leq \int^{M}_3 e^{\big[(1-\sigma_0)s^{\tau_0}\zeta_{M}(s)+\sqrt{\sigma_0(2-\sigma_0)s}\big]^2-s}ds 
  +\int^\infty_{3}e^{\big[(1-\sigma_0)C_{M}M^{1-2\tau_0}s^{\tau_0}+\sqrt{\sigma_0(2-\sigma_0) s}\big]^2-s}ds\\
  &\to \int^{M}_3 e^{s^{2\tau_0}\zeta_{M}(s)^2-s}ds 
  +\int^\infty_{3}e^{\big[C_{M}M^{1-2\tau_0}s^{\tau_0}\big]^2-s}ds\\
  &<+\infty\ \, \mbox{ as } \sigma_0\to 0^+.
  \end{align*}
  
  Therefore, if the desired result in the lemma does not hold, then we can find a positive sequence $\{\sigma_n\}$ decreasing to 0 as $n\to\infty$, and a positive sequence $\{s_n\}$ increasing to $\infty$ as $n\to\infty$
  such that 
  \begin{equation}\label{unbdd}
  \int^{s_n}_3 e^{A_n(s)^2 -s}ds 
  +\int^\infty_{s_n}e^{B_n(s)^2-s}ds\to\infty \mbox{ as $n\to\infty$},
  \end{equation}
  where
  \[\begin{cases}
  A_n(s) =(1\!-\!\sigma_n)s^{\tau_0}\zeta_{s_n}(s)\!+\!\sqrt{\sigma_n(2-\sigma_n)}\min\big\{ s^{\frac{1}2}\!, s^{\tau_0}s_n^{\frac 12\!-\!\tau_0}\!-\! s^{\frac{1}{2}} \big\},\\[2mm]
  B_n(s) =   (1\!-\!\sigma_n)C_{s_n}s_n^{1\!-\!2\tau_0}s^{\tau_0}\!+\!\sqrt{\sigma_n(2-\sigma_n)} (s^{\frac 12}\!-\!s^{\tau_0}s_n^{\frac 12 -\tau_0}).
  \end{cases}
  \]
  
  {\bf Step 1.} Estimate of $\displaystyle\int^{s_n}_3 e^{A_n(s)^2 -s}ds$.
  
  Fix $\delta\in(0,1)$ small. For $ 3\leq s\leq (1-\delta)s_n$, we have
  \begin{align*}
  A_n(s)&\leq (1-\sigma_n)s^{\tau_0}\zeta_{s_n}(s)\!+\!\sqrt{\sigma_n(2-\sigma_n)} s^{\frac{1}2}\\
  &= \Big[(1-\sigma_n)\sqrt{\frac {s^{1-2\tau_0}}{s_n^{1-2\tau_0}-2^{1-2\tau_0}}}+\sqrt{\sigma_n(2-\sigma_n)}\Big]s^{\frac 12}\\
  &\leq \Big[(1-\sigma_n)(1-\delta)^{\frac{1-2\tau_0}2}\sqrt{\frac {s_n^{1-2\tau_0}}{s_n^{1-2\tau_0}-2^{1-2\tau_0}}}+\sqrt{\sigma_n(2-\sigma_n)}\Big]s^{\frac 12}\\
  &=[(1-\delta)^{\frac{1-2\tau_0}2}+o(1)]s^{\frac 12}.
  \end{align*}
  It follows that
  \[
  A_n(s)^2 -s=[(1-\delta)^{1-2\tau_0}-1+o(1)]s\leq -\delta_0 s \mbox{ with $\delta_0=\frac{1-(1-\delta)^{1-2\tau_0}}2>0$ for all large $n$}.
  \]
  Therefore
  \[
  \int^{(1-\delta)s_n}_3 e^{A_n(s)^2 -s}ds\leq \int^{(1-\delta)s_n}_3 e^{-\delta_0 s}ds\leq e^{-3\delta_0}/\delta_0 \mbox{ for all large } n.
  \]
  
  Next we fix $r>0$ small so that 
  \[
  r\frac{1-2\tau_0}{\sqrt 2}<1.
  \]
  Then for
  $(1-\delta)s_n\leq s\leq (1-r\sqrt{\sigma_n})s_n$, we have
  \begin{align*}
  A_n(s)&\leq (1-\sigma_n)s^{\tau_0}\zeta_{s_n}(s)\!+\!\sqrt{\sigma_n(2-\sigma_n)} \Big[(1-\delta)^{\tau_0-\frac 12}-1\Big]s^{\frac{1}2}\\
  &= \Big\{(1-\sigma_n)\sqrt{\frac {s^{1-2\tau_0}}{s_n^{1-2\tau_0}-2^{1-2\tau_0}}}+\sqrt{\sigma_n(2-\sigma_n)}\Big[(1-\delta)^{\tau_0-\frac 12}-1\Big]\Big\}s^{\frac 12}\\ 
  &\leq\Big\{(1-\sigma_n)\Big(1-r\sqrt{\sigma_n}\Big)^{\frac{1-2\tau_0}2}\sqrt{\frac {s_n^{1-2\tau_0}}{s_n^{1-2\tau_0}-2^{1-2\tau_0}}}+\sqrt{\sigma_n(2-\sigma_n)}\Big[(1-\delta)^{\tau_0-\frac 12}-1\Big]\Big\}s^{\frac 12}\\
   &= \left\{1-\Big(\frac{1-2\tau_0}{2}r- \sqrt{2}\Big[(1-\delta)^{\tau_0-\frac 12}-1\Big]+o(1)\Big)\sqrt{\sigma_n}+[2^{-2\tau_0}+o(1)]s_n^{2\tau_0-1}\right\}s^{\frac 12}. 
  \end{align*}
  By choosing $\delta>0$ small enough, we have
  \[
  \sigma^*:=\frac{1-2\tau_0}{2}r- \sqrt{2}\Big[(1-\delta)^{\tau_0-\frac 12}-1\Big]>0,
  \]
  and hence, for $(1-\delta)s_n\leq s\leq (1-r\sqrt{\sigma_n})s_n$,
  \[
  A_n(s)^2-s\leq -\sigma^*\sqrt{\sigma_n} s+4^{1-\tau_0} s_n^{2\tau_0-1}s \leq -\sigma^*\sqrt{\sigma_n} s_n+4^{1-\tau_0} s_n^{2\tau_0}.
  \]
  It follows that
  \[
   \int_{(1-\delta)s_n}^{(1-r\sqrt{\sigma_n})s_n} e^{A_n(s)^2 -s}ds\leq \delta s_ne^{-\sigma^*\sqrt{\sigma_n} s_n+4^{1-\tau_0} s_n^{2\tau_0}} \to 0 \mbox{ as } n\to\infty,
  \]
  provided that 
  \begin{equation}\label{sigma_n}
  \sqrt{\sigma_n}s_n\geq s_n^\epsilon \mbox{ for some small $\epsilon>0$ and all large $n$.}
  \end{equation}
    
  For $(1-r\sqrt{\sigma_n})s_n\leq s\leq s_n$,
 we have
   \begin{align*}
  A_n(s)&\leq (1-\sigma_n)s^{\tau_0}\zeta_{s_n}(s)\!+\!\sqrt{\sigma_n(2-\sigma_n)} \Big[\Big(1-r\sqrt{\sigma_n}\Big)^{\tau_0-\frac 12}-1\Big]s^{\frac{1}2}\\
  &= (1-\sigma_n)C_{s_n}s^{1-\tau_0} +\sqrt{\sigma_n(2-\sigma_n)} \Big[\Big(1-r\sqrt{\sigma_n}\Big)^{\tau_0-\frac 12}-1\Big]s^{\frac{1}2}=:\tilde A_n(s),
  \end{align*}
  and so
\begin{align*}
\tilde A_n(s)
&\geq (1-\sigma_n)(1-r\sqrt{\sigma_n})^{\frac{1-2\tau_0}2}\frac{s_n^{\frac{1-2\tau_0}2}}{\sqrt{s_n^{1-2\tau_0}-2^{1-2\tau_0}}}s^{\frac 12} +o(1)s^{\frac 12}=[1+o(1)]s^{\frac 12}.
\end{align*}
\begin{align*}
\tilde A_n'(s)&=(1-\sigma_n)C_{s_n}(1-\tau_0)s^{-\tau_0} +\sqrt{\sigma_n(2-\sigma_n)} \Big[\Big(1-r\sqrt{\sigma_n}\Big)^{\tau_0-\frac 12}-1\Big]\frac 12s^{-\frac{1}2}\\
&\geq (1-\sigma_n)(1-\tau_0)(1-r\sqrt{\sigma_n})^{-\tau_0}\frac{s_n^{\frac{1-2\tau_0}2}}{\sqrt{s_n^{1-2\tau_0}-2^{1-2\tau_0}}}s^{-\frac 12} +o(1)s^{-\frac 12}\\
&=[1-\tau_0+o(1)]s^{-\frac 12},
\end{align*}
\begin{align*}
\left(e^{\tilde A_n(s)^2-s}\right)'&=e^{\tilde A_n(s)^2-s}[2\tilde A_n(s)\tilde A'_n(s)-1]\\[2mm]
&  \geq\  [1-2\tau_0+o(1)]e^{\tilde A_n(s)^2-s} 
\  \geq e^{\tilde A_n(s)^2-s}\ \,\mbox{ for all large } n.
\end{align*}
We thus have
\[\begin{aligned}
\int^{s_n}_{(1-r\sqrt{\sigma_n})s_n} e^{A_n(s)^2 -s}ds&\leq\int^{s_n}_{(1-r\sqrt{\sigma_n})s_n} e^{\tilde A_n(s)^2 -s}ds\\[2mm]
&\leq\int^{s_n}_{(1-r\sqrt{\sigma_n})s_n}\left(e^{\tilde A_n(s)^2-s}\right)'ds
\leq e^{\tilde A_n(s_n)^2-s_n}.
\end{aligned}
\]
Since
\begin{align*}
\tilde A_n(s_n)&\leq  (1-\sigma_n)s_n^{\tau_0}\zeta_{s_n}(s_n)\!+\!\sqrt{\sigma_n(2-\sigma_n)} \Big[\Big(1-r\sqrt{\sigma_n}\Big)^{\tau_0-\frac 12}-1\Big]s_n^{\frac{1}2}\\
&=(1-\sigma_n)\frac{s_n^{\frac{1-2\tau_0}2}}{\sqrt{s_n^{1-2\tau_0}-2^{1-2\tau_0}}}s_n^{\frac 12}\!+\!\sqrt{\sigma_n(2-\sigma_n)} \Big[\Big(1-r\sqrt{\sigma_n}\Big)^{\tau_0-\frac 12}-1\Big]s_n^{\frac{1}2}\\
&=\left\{  (1-\sigma_n)\Big[1+\Big(\frac 12+o(1)\Big)2^{1-2\tau_0}s_n^{2\tau_0-1} \Big]+ \Big(\frac{1-2\tau_0}{\sqrt 2} r+o(1)\Big)\sigma_n \right\}s_n^{\frac{1}2}\\
&=\left\{  1-\sigma_n\Big[1-\frac{1-2\tau_0}{\sqrt {2}} r+o(1)\Big]+\Big[2^{-2\tau_0}+o(1)\Big]s_n^{2\tau_0-1} \right\}s_n^{\frac 12}\\
&\leq \left( 1+2^{1-2\tau_0}s_n^{2\tau_0-1} \right)s_n^{\frac 12}\ \mbox{ for all large $n$ due to the choice of } r,
\end{align*}
we obtain
\[
\tilde A_n(s_n)^2-s_n\leq 2^{3-2\tau_0} s_n^{2\tau_0} \mbox{ for all large } n,
\]
and
\[\begin{aligned}
\int^{s_n}_{(1-r\sqrt{\sigma_n})s_n} e^{A_n(s)^2 -s}ds
\leq e^{\tilde A_n(s_n)^2-s_n}\leq e^{2^{3-2\tau_0} s_n^{2\tau_0}}\to 1\ \,  \mbox{ as } n\to\infty.
\end{aligned}
\]
\medskip

If \eqref{sigma_n} does not hold, then for any small $\epsilon>0$, there is a subsequence of $\{\sigma_n\}$, still denoted by itself for simplicity of notation, such that
\[
\sqrt{\sigma_n}\leq s_n^{\epsilon-1} \mbox{ for all } n\geq 1.
\]
Therefore, for fixed $\sigma\in (1-\epsilon,1)$ and $s\in [(1-\delta)s_n, (1-s_n^{-\sigma})s_n]$,
\begin{align*}
A_n(s)&= s^{\tau_0}\zeta_{s_n}(s)+O(s_n^{\epsilon-1}) s^{\frac{1}{2}}\\
&\leq (1-s_n^{-\sigma})^{\frac{1-2\tau_0}2}\frac{s_n^{\frac{1-2\tau_0}2}}{\sqrt{s_n^{1-2\tau_0}-2^{1-2\tau_0}}}s^{\frac 12}+O(s_n^{\epsilon-1}) s^{\frac{1}{2}}\\
&=\big[ 1-{\frac{1-2\tau_0}2}s_n^{-\sigma}+o(s_n^{-\sigma})\big]s^{\frac 12}.
\end{align*}
Hence
\[
A_n(s)^2-s\leq \big[-(1-2\tau_0)s_n^{-\sigma}+o(s_n^{-\sigma})\big]s\leq -s_n^{-\sigma}s\leq -(1-\delta)s_n^{1-\sigma} \mbox{ for all large } n.
\]
It follows that
\[
\int_{(1-\delta)s_n}^{(1-s_n^{-\sigma})s_n}e^{A_n(s)^2-s}ds\leq \delta s_n e^{-(1-\delta)s_n^{1-\sigma}} \to 0 \mbox{ as } n\to\infty.
\]

For  $s\in [(1-s_n^{-\sigma})s_n, s_n]$, we have
\[
A_n(s)\leq (1-\sigma_n)C_{s_n}s^{1-\tau_0} +\sqrt{\sigma_n(2-\sigma_n)} \Big[\Big(1-s_n^{-\sigma}\Big)^{\tau_0-\frac 12}-1\Big]s^{\frac{1}2}=: A^*_n(s),
\]
and
\[
A^*_n(s)=[1+o(1)]s^{1/2},\ A^*_n(s)'=[1-\tau_0+o(1)]s^{-1/2}.
\]
Therefore
\[
\left(e^{A^*_n(s)^2-s}\right)'=e^{A^*_n(s)^2-s}[2A^*_n(s)A^*_n(s)'-1]=e^{A_n(s)^2-s}[1-2\tau_0+o(1)]\geq e^{A_n(s)^2-s} 
\]
 for all large $n$. It follows that
\[
\int_{(1-s_n^{-\sigma})s_n}^{s_n} e^{A_n(s)^2-s}ds\leq \int_{(1-s_n^{-\sigma})s_n}^{s_n} e^{A^*_n(s)^2-s}\leq \int_{(1-s_n^{-\sigma})s_n}^{s_n} \left(e^{A^*_n(s)^2-s}\right)'ds\leq e^{A^*_n(s_n)^2-s_n}
\]
for all large $n$.
We now calculate
\begin{align*}
A^*_n(s_n)&=(1-\sigma_n)C_{s_n}s_n^{1-\tau_0} +\sqrt{\sigma_n(2-\sigma_n)} \Big[\Big(1-s_n^{-\sigma}\Big)^{\tau_0-\frac 12}-1\Big]s_n^{\frac{1}2}\\
&=\Big\{(1-\sigma_n)\frac{s_n^{\frac{1-2\tau_0}2}}{\sqrt{s_n^{1-2\tau_0}-2^{1-2\tau_0}}}+\big[\frac{2\tau_0-1}{\sqrt 2}+o(1)\big]\sqrt{\sigma_n} s_n^{-\sigma}\Big\}s_n^{1/2}\\
&=\Big\{1+O(s_n^{2\epsilon-2})+O(s_n^{2\tau_0-1})+O(s_n^{\epsilon-1-\sigma})\Big\}s_n^{1/2}.
\end{align*}
Therefore
\[
A^*_n(s_n)^2-s_n=\Big\{O(s_n^{2\epsilon-2})+O(s_n^{2\tau_0-1})+O(s_n^{\epsilon-1-\sigma})\Big\}s_n=o(1)
\]
provided that $\epsilon>0$ is sufficiently small, which gives
\[
\int_{(1-s_n^{-\sigma})s_n}^{s_n} e^{A_n(s)^2-s}ds\leq e^{A^*_n(s_n)^2-s_n}=e^{o(1)}\to 1 \mbox{ as } n\to\infty.
\]

Summarising, we have proved
\[
\liminf_{n\to\infty}\int_3^{s_n} e^{A_n(s)^2-s}ds<\infty.
\]

 {\bf Step 2.} Estimate of
$\displaystyle
\int^\infty_{s_n} e^{B_n(s)^2-s}ds.
$

Clearly, for $s\geq s_n$,
\[B_n(s)\leq  (1-\sigma_n)C_{s_n}s_n^{1-2\tau_0} s^{\tau_0}+\sqrt{2\sigma_n}s_n^{\tau_0}(s^{\frac 12 -\tau_0}-s_n^{\frac 12-\tau_0}) =: \tilde B_n(s).
\]
Fix $\epsilon_0>0$ small. For $s_n\leq s\leq (1+\epsilon_0)s_n$, we have
\begin{align*}
\tilde B_n(s)\geq (1-\sigma_n)(1+\epsilon_0)^{\tau_0}\frac{s_n^{\frac{1-2\tau_0}2}}{\sqrt{s_n^{1-2\tau_0}-2^{1-2\tau_0}}}s_n^{1/2}=[(1+\epsilon_0)^{\tau_0}+o(1)]s_n^{1/2},
\end{align*}
and
\begin{align*}
-\tilde B_n'(s)&=-\tau_0 (1-\sigma_n)\frac{s_n^{1-2\tau_0}}{\sqrt{s_n^{1-2\tau_0}-2^{1-2\tau_0}}}s^{\tau_0-1}+\sqrt{2\sigma_n}s_n^{\tau_0}s^{-\frac 12-\tau_0}(\frac 12-\tau_0)\\
&\geq -\tau_0(1-\sigma_n)(1+\epsilon_0)^{\tau_0-1}\frac{s_n^{\frac{1-2\tau_0}2}}{\sqrt{s_n^{1-2\tau_0}-2^{1-2\tau_0}}}s_n^{-1/2}+o(1)s_n^{-1/2}\\
&=[-\tau_0(1+\epsilon_0)^{\tau_0-1}+o(1)]s_n^{-1/2}.
\end{align*}
It follows that
\[
-\tilde B_n(s)\tilde B_n'(s)\geq -\tau_0(1+\epsilon_0)^{2\tau_0-1}+o(1)\geq -\tau_0(1+2\epsilon_0)^{2\tau_0-1} \mbox{ for all large } n,
\]
and therefore
\[
\left(e^{\tilde B_n(s)^2-s}\right)'=e^{\tilde B_n(s)^2-s}[2\tilde B_n(s)\tilde B_n'(s)-1]\leq -\xi_0 e^{\tilde B_n(s)^2-s}
\]
with
\[ \xi_0:=-2\tau_0(1+2\epsilon_0)^{2\tau_0-1}+1>0.
\]
We thus obtain
\[
\int^{(1+\epsilon_0)s_n}_{s_n} e^{B_n(s)^2-s}ds\leq \int^{(1+\epsilon_0)s_n}_{s_n} e^{\tilde B_n(s)^2-s}ds\leq -\xi_0^{-1}\int^{(1+\epsilon_0)s_n}_{s_n} \left(e^{\tilde B_n(s)^2-s}\right)'ds\leq \xi_0^{-1}e^{\tilde B_n(s_n)^2-s_n}
\]
for all large $n$. Clearly
\[
\tilde B_n(s_n)=(1-\sigma_n)C_{s_n}s_n^{1-2\tau_0} s_n^{\tau_0}\leq\frac{s_n^{\frac{1-2\tau_0}2}}{\sqrt{s_n^{1-2\tau_0}-2^{1-2\tau_0}}}s_n^{1/2}=\big[1+O(s_n^{2\tau_0-1})\big]s_n^{1/2}.
\]
Therefore
\[
\tilde B_n(s_n)^2-s_n=O(s_n^{2\tau_0-1}) s_n=o(1),
\]
and thus
\[
\int^{(1+\epsilon_0)s_n}_{s_n} e^{B_n(s)^2-s}ds\leq\xi_0^{-1}e^{\tilde B_n(s_n)^2-s_n}
\leq \xi_0^{-1} e^{o(1)}\to \xi_0^{-1} \mbox{ as } n\to\infty.
\]

For $s\geq (1+\epsilon_0)s_n$, we have
\begin{align*}
B_n(s)&\leq (1-\sigma_n)C_{s_n}s_n^{1-2\tau_0} s^{\tau_0}+\sqrt{2\sigma_n}s^{\frac 12}\\[1mm]
&\leq (1-\sigma_n)(1+\epsilon_0)^{\tau_0}C_{s_n}s_n^{1-2\tau_0} s_n^{\tau_0}+\sqrt{2\sigma_n}s^{\frac 12}\\
&=\big[(1+\epsilon_0)^{\tau_0}+o(1)\big]\frac{s_n^{\frac{1-2\tau_0}2}}{\sqrt{s_n^{1-2\tau_0}-2^{1-2\tau_0}}}s_n^{1/2}+o(1)s^{1/2}\\
&=\big[(1+\epsilon_0)^{\tau_0}+o(1)\big]s_n^{1/2}+o(1)s^{1/2}.
\end{align*}
It follows that
\[
B_n(s)^2-s\leq \big[(1+\epsilon_0)^{2\tau_0}+o(1)\big]s_n-\big[1+o(1)\big]s\leq (1+\epsilon_0/2)^{2\tau_0}(s_n-s) \mbox{ for all large } n.
\]
Therefore,
\begin{align*}
\int_{(1+\epsilon_0)s_n}^\infty e^{B_n(s)^2-s}ds&\leq\int_{(1+\epsilon_0)s_n}^\infty e^{(1+\epsilon_0/2)^{2\tau_0}(s_n-s)}ds
\\[2mm]&=(1+\epsilon_0/2)^{-2\tau_0}e^{-(1+\epsilon_0/2)^{2\tau_0}\epsilon_0 s_n}\to 0\quad{\rm as}\ \, n\to\infty.
\end{align*}

We thus obtain
\[
\limsup_{n\to\infty}\int_{s_n}^\infty e^{B_n(s)^2-s}ds\leq \xi_0^{-1}.
\]
Combining this with the conclusion proved in Step 2a we obtain
\[
\liminf_{n\to\infty}\left[\int_3^{s_n} e^{A_n(s)^2-s}ds+\int_{s_n}^\infty e^{B_n(s)^2-s}ds\right]<\infty.
\]
But this clearly contradicts \eqref{unbdd}. The proof of the lemma is now complete.
\hfill$\Box$

  \begin{lemma}\label{pr 40}
  Let $\mu>-\frac14$.  Then
    there exist a sequence of radially symmetric functions $u_n\in \cH^1_{\mu, 0}(B_1)$  such that 
$$\|u_n\|_{\mu}^2=\int_{B_1}\Big(|\nabla u_n(x)|^2+\mu \frac{u_n(x)^2 }{|x|^2 (-\ln |x| )^2} \Big) dx= 1  $$
and 
$$\lim_{n\to\infty} \int_{B_1} e^{\alpha |u_n|^2}dx\to\infty\ \mbox{ for  any
  $\alpha>m_\mu =4\pi    \sqrt{1+4\mu}$}. $$
 \end{lemma}
\noindent{\bf Proof. } It is convenient to do the construction for a transformed form of $u_n$.
Let   
$$u_n(r)=(-\ln r)^{\tau_0} v_n(r)\quad{\rm for}\ \, r\in(0,1);$$
then $v_n(0)=0$ and $v_n(r)$ also vanishes in the small neighbourhood of $r=1$. 
 Moreover,
\begin{eqnarray*} 
 \|u_n\|^2_\mu&=&2\pi \int^1_0 \Big(u_n'(r)^2+\mu \frac{u_n(r)^2}{r^2(-\ln r)^2}\Big)  rdr
 \\[2mm]&=&  2\pi \int^1_0  v_n'(r)^2  (-\ln r)^{2\tau_0} rdr. 
    \end{eqnarray*} 
 Let
  $$  r=e^{-\frac12 t}  \quad {\rm and}\quad v_n(r)=w_n(t);$$
 then
 \begin{eqnarray*} 
  \|u_n\|^2_\mu   =  2^{2-2\tau_0}\pi \int^\infty_0   w_n'(t)^2 t^{2\tau_0} dt=:\cJ(w_n)
    \end{eqnarray*}
 and
  $$
 \int_{B_1} e^{\alpha |u_n|^2}dx=2\pi \int^1_0  e^{\alpha u_n^2}r dr = \pi   \int^\infty_0  \exp( \alpha  2^{-2\tau_0} t^{2\tau_0}w_n^2-t )  dt=:\cN(w_n).$$

We now construct  functions $\{w_n\}$ such that  $\cJ(w_n)=1$ and $\cN(w_n)\to\infty$.   Let
 \begin{eqnarray}  \label{test fd}
   \zeta_\sigma(t): =\left\{ \arraycolsep=1pt
\begin{array}{lll}
t^\sigma\ \  &{\rm for}\ \,   t\in[0,1),\\[2mm]
1\ \  &{\rm for}\ \,   t\in[1  +\infty),
\end{array}
\right.  
\end{eqnarray} 
with $\sigma>\frac12-\tau_0$ to be determined later. 
 Then define, for all large integer $n$ and some constant $\nu>0$ to be specified below,
$$w_n(t):=\nu n^{\frac12-\tau_0}\zeta_\sigma(n^{-1}t) . $$
 Direct computation gives
  \begin{eqnarray*} 
   \int^\infty_0   w_n'(t)^2 t^{2\tau_0}    dt 
   &=&   \nu^2 \sigma^2 \int^n_0  n^{1-2\tau_0-2\sigma}   t^{2\tau_0+2\sigma-2} dt 
  \\[2mm]&=&     \frac{\nu^2 \sigma^2 }{2\tau_0+2\sigma-1}.  
    \end{eqnarray*}
 Therefore
   \begin{eqnarray} \label{est J}
 \cJ(w_n)
   =  2^{2-2\tau_0}\pi \frac{\nu^2 \sigma^2 }{2\tau_0+2\sigma-1}   = 1
    \end{eqnarray}
provided that
\[
\nu=    \big(2^{2-2\tau_0}\pi f_\mu(\sigma) \big)^{-1/2},
\]
where
\[
f_\mu(\sigma):= \frac{  \sigma^2  }{2\sigma+2\tau_0-1} \quad \mbox{for } \sigma>\frac12-\tau_0.
\]

We next estimate $\cN(w_n)$. 
To this end,  we first choose an optimal $\sigma$. It is easily seen that $f_\mu(\sigma)$
 achieves  its minimum   at 
 \begin{eqnarray*} 
\sigma= \sigma_\mu =:   1-2\tau_0=\sqrt{1+4\mu}  
    \end{eqnarray*} 
 with     
 $$f_\mu (\sigma_\mu)=1-2\tau_0=\sqrt{1+4\mu}.  $$
  Then (\ref{est J}) holds with $\sigma=\sigma_\mu$  and
 \begin{eqnarray*} 
 \nu  =\nu_\mu:=   \big(2^{2-2\tau_0}\pi f_\mu(\sigma_\mu) \big)^{-1/2}. 
 \end{eqnarray*} 
 
With $(\sigma,\nu)=(\sigma_\mu, \nu_\mu)$ in the definition of $w_n$, by   direct computation, we have
 \begin{eqnarray*} 
\cN(w_n)
 &>&  \pi   \int^n_0  \exp( \alpha  \nu_\mu^2 2^{-2\tau_0} n^{1-2\tau_0-2\sigma_\mu} t^{2\sigma_\mu+2\tau_0} -t )  dt
 \\[2mm]&>&  \pi  \int^{n}_{(1-\epsilon_0)n} \exp( \alpha \nu_\mu^2 2^{-2\tau_0} n^{1-2\tau_0-2\sigma_\mu} (1-\epsilon_0)^{2\sigma_\mu+2\tau_0}n^{2\sigma_\mu+2\tau_0} -n )  dt 
  \\[2mm]&=&  \pi e^{( \alpha \nu_\mu^2 2^{-2\tau_0}   (1-\epsilon_0)^{2\sigma_\mu+2\tau_0} -1)n }   \int^{n}_{(1-\epsilon_0)n}  dt  
  \\[2mm]&\geq & \pi \epsilon_0 n
 \to   \infty \quad{\rm as}\ \, n\to\infty
    \end{eqnarray*}  
if $$\alpha \nu_\mu^2 2^{-2\tau_0}   (1-\epsilon_0)^{2\sigma_\mu+2\tau_0} -1\geq0,$$
   which holds when 
      \begin{eqnarray} \label{est a}
  \alpha >  2^{2\tau_0} \nu_\mu^{-2} = 4\pi  \sqrt{1+4\mu} \mbox{ and } 0<\epsilon_0\ll 1.  
 \end{eqnarray}   
  \hfill$\Box$ \medskip
  
 Theorem \ref{teo 1} clearly follows directly from Lemmas \ref{teo 3.1}  and \ref{pr 40}.

    \subsection{Proof of Theorem \ref{teo 2} }

For the case $\mu=-\frac14$,  if we do the transform that  $r=e^{-\frac12 t}$ and $w_0(t)=v_0(r)=(-\ln r)^{-\tau_0}u_0(r)$,  it fails to provide a useful  estimate of the form
(\ref{eqvi 1-2.1}), since $1-2\tau_0= 0$ in this case.    Instead, we will find a new transformation.

   \begin{lemma}\label{teo 3.2}
Let   $u_0\in \cH^1_{{\rm rad},-1/4, 0}(B_1)$ satisfy
$$\|u_0\|_{-1/4}^2=\int_{B_1}\Big(|\nabla u_0(x)|^2-\frac{u_0(x)^2}{4|x|^2 (-\ln |x| )^2} \Big) dx\leq 1.  $$
Then for $p\in(0,1)$ and $\alpha>0$,     there exists a constant $C_{p,\alpha}>0$ independent of $u_0$ such that 
$$\int_{B_1} e^{\alpha |u_0|^{p }}dx\leq  C_{p,\alpha}.$$  
 \end{lemma}
\noindent {\bf Proof. }   Replacing $u_0$ by $|u_0|$, we may assume that $u_0\in \cH^1_{{\rm rad},\mu, 0}(B_1)$ is nonnegative.  By the embedding theorem we see
 that  $u_0$ is continuous in $[0,1]$.  
Moreover, since the space  $C_{{\rm rad}, c}^\infty(B_1)$ is dense in  $\cH^1_{{\rm rad}, \mu, 0}(B_1)$, we may further assume that $u_0\in C_{{\rm rad}, c}^\infty(B_1)$. 
By Theorem \ref{embedding}, 
$$\|u_0\|_{L^2(B_1)}\leq C\|u_0\|_{\mu}\leq C,$$
 with $C$ independent of $u_0$. Fix $r_1< r_2$ such that
  $0<r_1\leq 1/4,\  3/4\leq r_2<1$. Then
  \[
  C\geq \|u_0\|_{L^2(B_1)}\geq \|u_0\|_{L^2(B_{r_2}\setminus B_{r_1})}\geq |B_{r_2}\setminus B_{r_1})|^{1/2}\inf_{B_{r_1}\setminus B_{r_2} } u_0.
  \]
 Therefore there exists $r_0\in [r_1, r_2]$ such that
 \[
 u_0(r_0)\leq \tilde C:=C|B_{r_2}\setminus B_{r_1})|^{-1/2}.
 \]

    Let 
 $$ u_0(r)=(-\ln r)^{\frac12} v_0(r)\quad{\rm for}\ \, r\in(0,1),$$
  $$   v_0(r)=w_0(t),\  r=e^{\frac12(1-  e^t) },\ r_0=e^{\frac12(1-  e^{t_0}) }.$$
 Then $dr=-\frac12e^{\frac12(1-  e^t) } e^tdt,$
 $$w_0(0)=\lim_{t\to0^+} w_0(t)=\lim_{r\to1^-}v_0(r)=0,\qquad w_0(t_0)\leq \frac{2\tilde C}{e^{t_0}-1} ,$$
and  direct computation shows that
 \begin{equation}\label{I(u_0)}\begin{aligned}
  1\geq \|u_0\|_{-1/4}^2    &=    2\pi \int^1_0   v_0'(r)^2 r(-\ln r)  dr
  \\[2mm]&= 2\pi \int^{+\infty}_0  w_0'(t)^2 \,  \frac{e^t-1}{e^t}   dt 
  \\[2mm]& \geq 2\pi \big(1-e^{-t_0}\big)\int^{+\infty}_{t_0}  w_0'(t)^2 \,      dt+2\pi e^{-t_0}\int^{t_0}_0  w_0'(t)^2tdt,
  \end{aligned}
    \end{equation}
     since $1-e^{-t}\geq e^{-t_0}t $ for $t\in(0, t_0)$.
  Moreover, 
  \begin{eqnarray} \label{chac}  
   2\pi \int^1_0  e^{\alpha |u_0|^p } r dr
 =   \pi   \int^{+\infty}_0  \Phi_0(t)  dt,
   \end{eqnarray}
  where 
  $$\Phi_0(t):=e^{\alpha \frac{(e^t-1)^p}{2^p} |w_0|^p+1-e^t +t }. $$
   
   By \eqref{I(u_0)} we have
   $$\int^{+\infty}_{t_0}  w_0'(t)^2 \,      dt\leq b_0^2:= \big[2\pi \big(1-e^{-t_0}\big)  \big]^{-1}.$$
 It follows that, for $t>t_0$,
   \begin{eqnarray*} 
  w_0(t)    \leq    \int^t_{t_0}   w_0'(s)  ds+w_0(t_0)&\leq& \Big(\int^t_{t_0}  w_0'(s)^2     ds \Big)^{\frac12} \Big(\int^t_1      ds\Big)^{\frac12}+w_0(t_0)
  \\[2mm]&\leq & b_0   t^{\frac12}+\tilde C,  
    \end{eqnarray*}
  and hence, for $p\in (0,1)$
 \begin{eqnarray*} 
   \pi   \int^{+\infty}_{t_0}   \Phi_0(t)  dt
 &=&  \pi   \int^{+\infty}_{t_0}  e^{\alpha \frac{(e^t-1)^p}{2^p} |w_0|^p+1-e^t +t }  dt
 \\[2mm]&\leq &  \pi  \int^{+\infty}_{t_0}  e^{\alpha 2^{-p}  e^{pt} \big[  b_0   t^{\frac12}+\tilde C\big]^p +1-e^t +t }  dt\\
 & \leq &  C_{p,\alpha}.
    \end{eqnarray*} 
    where we have used the fact that $t_0$  can be bounded (from above and below) by positive numbers depending on $r_1, r_2$ only, via $r_1\leq r_0\leq r_2$.

 By \eqref{I(u_0)} we also have
   $$\int^{t_0}_0  w_0'(t)^2 t\,      dt\leq b_1^2:=\big[{2\pi}(1-e^{-t_0})  \big]^{-1}. $$  
 Thus, for $t\in (0, t_0)$, 
   \begin{eqnarray*} 
 w_0(t)   &\leq&  w_0(t_0) - \int_t^{t_0}   w_0'(s)  ds 
 \\[2mm]&\leq& w_0(t_0)+ \Big(\int_t^{t_0}  w_0'(s)^2  s   ds \Big)^{\frac12} \Big(\int_t^{t_0}  s^{-1}   ds\Big)^{\frac12} 
      \leq   \tilde C+  b_1 \sqrt{\ln \frac{t_0}t}. 
    \end{eqnarray*}
Since $p\in(0,1)$, for any $\epsilon>0$, there exists $t_1\in(0,t_0)$ such that
$$\Big( \tilde C+  b_1 \sqrt{\ln \frac{t_0}t}\Big)^{\frac{p}2}\leq \epsilon \ln\frac{t_0}t\quad{\rm for}\ \, t\in(0,t_1). $$
Now we take $\epsilon=\frac12 \frac{1}{\alpha 2^{-p}  e^{pt_0}   }$ and obtain
 \begin{eqnarray*} 
   \pi   \int^{t_0}_0   \Phi_0(t)  dt
  &\leq &  \pi  \int^{t_0}_0  e^{\alpha 2^{-p}  e^{pt} \big[  b_1  \sqrt{\ln\frac{t_0}t} +\tilde C\big]^p +1-e^t +t }  dt
  \\[2mm] &\leq & \pi  \int^{t_1}_0    e^{ \frac12  \ln\frac{t_0}{t  } } dt+ \pi  \int^{t_0}_{t_1}  e^{\alpha 2^{-p}  e^{pt} \big[  b_1  \sqrt{\ln\frac{t_0}t} +\tilde C\big]^p +1-e^t +t }  dt\\[2mm]
  &=& \pi \int^{t_1}_0  \frac{\sqrt{t_0}}{\sqrt t}  dt+ \pi  \int^{t_0}_{t_1}  e^{\alpha 2^{-p}  e^{pt} \big[  b_1  \sqrt{\ln\frac{t_0}t} +\tilde C\big]^p +1-e^t +t }  dt\\
 &\leq &  C_{p,\alpha},
    \end{eqnarray*} 
    where we have used the fact that $t_0$ and $t_1$ can be bounded (from above and below) by positive numbers depending on $r_1, r_2$ and $p$ only. 
Therefore  
 \begin{eqnarray*} 
  2\pi \int^1_0  e^{\alpha |u_0|^p } r dr
  =   \pi   \int^{+\infty}_0  \Phi_0(t)   dt \leq  2 C_{p,\alpha},
    \end{eqnarray*} 
 and the proof is complete.    \hfill$\Box$ \medskip
  
  \begin{lemma}\label{pr 41}
  There exist a sequence of functions $u_n\in \cH^1_{{\rm rad}, -1/4, 0}(B_1)$  such that 
$$\int_{B_1}\Big(|\nabla u_n(x)|^2- \frac{u_n(x)^2 }{4|x|^2 (-\ln |x| )^2} \Big) dx\leq 1  $$
and  
$$\lim_{n\to\infty}\int_{B_1} e^{\alpha |u_n|^p}dx=\infty\quad\mbox{ for any } \alpha>0,\ p\geq 1.$$
 \end{lemma}
\noindent{\bf Proof. } We will construct a sequence of radially symmetric functions  $\{u_n\}\subset C_c^\infty(B_1)$, which stays inside the unit ball of $\cH^1_{\mu_0, 0}(B_1)$, but  $\int_{B_1} e^{\alpha |u_n|^p}dx\to\infty$
for any $p\geq 1$ and any $\alpha>0$.

Using  the same transformations as in the proof of Lemma \ref{teo 3.2},  we have
  \begin{eqnarray*} 
  \|u_n\|_{-1/4}^2 = 2\pi \int^{+\infty}_0  w_n'(t)^2 \,  \frac{e^t-1}{e^t}   dt
    \end{eqnarray*}
    and
  \begin{eqnarray*} \int_{B_1} e^{\alpha |u_n|^p}dx=  \pi   \int^{+\infty}_0  e^{\alpha \big[ \frac{ e^t-1 }{2 }\big]^p |w_n |^p+1-e^t +t }  dt.   
  \end{eqnarray*}
 
 Given $\kappa>4$, we let $\eta_\kappa:(-\infty, +\infty)\to[0,1]$ be a  smooth function such that 
\begin{eqnarray} \label{cutoff f1} 
  \eta_\kappa(t): =\left\{ \arraycolsep=1pt
\begin{array}{lll}
 0\ \ \ \ &{\rm for}\ \,   t\in(-\infty,1),
\\[2mm]
1\ \  &{\rm for}\ \,   t\in(2,\kappa-1),\\[2mm]
0\ \  &{\rm for}\ \,   t\in(\kappa,+\infty).
\end{array}
\right.  
    \end{eqnarray} 
 Denote
   $$b_1(\kappa):=\int^{\kappa}_1  \eta_\kappa (s) ds\quad{\rm and}\quad b_2(\kappa):=\sqrt{4\pi \int^{\kappa}_1  \eta_\kappa(s)^2 ds};$$
then 
\begin{eqnarray*}  
 \frac{b_1(\kappa)}{b_2(\kappa)} &\geq &    \frac{\int_2^{\kappa-1} 1ds }{\sqrt{4\pi\int_1^{\kappa} 1ds} }
 =   \frac{ \kappa-3  }{\sqrt{4\pi(\kappa-1)} } \to+\infty\quad{\rm as}\ \, \kappa\to+\infty.
    \end{eqnarray*} 
Define
\begin{eqnarray} \label{cutoff f2}  
\phi_\kappa(t):=  \eta_\kappa(t)-\eta_\kappa( t-2\kappa) \quad{\rm for}\ \, t\in \R,
  \end{eqnarray} 
 and
  $$w_\kappa(t):=\frac1{b_2(\kappa)} \int_0^t\phi_\kappa (s)ds\quad{\rm for}\ \, t\geq0.$$
  It is easily seen that $w_\kappa$ is nonnegative and $w_\kappa\in C_c^\infty((0,\infty))$.  
 Moreover, 
 \begin{eqnarray*}  
 2\pi \int^{+\infty}_0  w_\kappa '(t)^2 \,  \frac{e^t-1}{e^t}   dt &=&  \frac{2\pi}{b_2(\kappa)^2}   \int^{+\infty}_0  \phi_\kappa (t)^2 \,  \frac{e^t-1}{e^t}   dt 
 \\[2mm] &\leq &  \frac{2\pi}{b_2(\kappa)^2}   \int^{+\infty}_0  \phi_\kappa(t)^2 dt
  =  1
  \end{eqnarray*} 
 and  for $ t\in[\kappa, 2\kappa]$, 
 $$w_\kappa(t)=\frac1{b_2(\kappa)}\int_0^{\kappa} \phi_\kappa (s)ds=\frac{b_1(\kappa)}{b_2(\kappa)},$$
 \begin{eqnarray*} 
   \pi   \int^{+\infty}_0  e^{\alpha \big[\frac{ e^t-1 }{2 }\big]^p |w_\kappa |^p+1-e^t +t }  dt&>& \pi   \int^{2\kappa}_{\kappa}   e^{\alpha \big[\frac{ e^t-1 }{2 }\big]^p |w_\kappa|^p+1-e^t +t }  dt
\\[2mm]  &= &  \pi   \int^{2\kappa}_{\kappa}   e^{\alpha \big[\frac{ e^t-1 }{2 }\big]^p \big[\frac{b_1(\kappa)}{b_2(\kappa)}\big]^p+1-e^t +t }  dt
 \\[2mm]&\to & +\infty\quad{\rm as}\ \ \kappa\to+\infty, 
     \end{eqnarray*}   
 which ends the proof.  \hfill$\Box$\medskip

Clearly Theorem \ref{teo 2}
 follows from Lemmas \ref{teo 3.2}  and  \ref{pr 41}  directly.

  \setcounter{equation}{0}
\section{Trudinger-Moser type inequalities for general functions}

\subsection{The case $\mu>-\frac14$ }

\begin{lemma}\label{lm super c}
Let   $\mu\in (-\frac14, 0)$,  $V:(0,1)\to[0,+\infty)$ be a continuous  function 
verifying  \eqref{con V-l}. Then there exists a sequence $\{u_n\}\subset \cH^1_{V, \mu, 0}(B_1)$
such that $\|u_n\|_{V,\mu}= 1$ and
  for any $\alpha>m_\mu$, $r_0\in (0,1)$,
$$\lim_{n\to\infty}\int_{B_{r_0}} e^{\alpha |u_n|^2}dx=\infty.$$
\end{lemma}
\noindent{\bf Proof.} As before, for $u_n\in C^{0,1}_{\rm rad, c}(B_1)$, under the transformation 
$$u_n(r)=(-\ln r)^{\tau_0} v_n(r)\quad{\rm for}\ \, r\in(0,1),$$
we have
\begin{eqnarray*} 
 \|u_n\|^2_{V,\mu} &=&2\pi \int^1_0 \Big(u_n'(r)^2+\mu Vu_n(r)^2  \Big)  rdr
 \\[2mm]&=&  2\pi \int^1_0  v_n'(r)^2  (-\ln r)^{2\tau_0} rdr + 2\pi\mu \int^1_0 \Big(V(r)-\frac{1}{r^2(-\ln r)^2}\Big) v_n(r)^2(-\ln r)^{2\tau_0}   rdr. 
    \end{eqnarray*} 
    
By (\ref{con V-l}), given small $\epsilon>0$, there exists $r_\epsilon\in(0,1)$ such that 
 \begin{eqnarray}\label{con V-l-c}
  |V(r) {r^2(-\ln r)^2}-1|\leq \epsilon\quad {\rm for}\ r\in(0,r_\epsilon].  
 \end{eqnarray}
   
 Let
  $$  r_\epsilon=e^{-\frac12 (t+t_\epsilon)}  \quad {\rm and}\quad v_n(r)=w_n(t),$$
where $t_\epsilon=-2\ln r_\epsilon$ so that $w_n(-t_\epsilon)=v_n(1)=0$.
Then
 $$
 2\pi \int^{r_0}_0  e^{\alpha u_n^2}r dr = \pi   \int^\infty_{-t_\epsilon-2\ln r_0}  e^{ \alpha  2^{-2\tau_0} (t+t_\epsilon)^{2\tau_0}w_n^2-(t+t_\epsilon) }  dt=:\tilde\cN(w_n)$$
 and
 \begin{eqnarray*} 
  2\pi \int^1_0  u_n'(r)^2   rdr=2\pi \int^1_0  v_n'(r)^2  (-\ln r)^{2\tau_0} rdr    = 2^{2-2\tau_0}\pi\int^\infty_ {-t_\epsilon}  w_n'(t)^2(t+t_\epsilon)^{2\tau_0} dt.
    \end{eqnarray*}

Next we construct the functions $\{w_n\}$.   Let
 \begin{eqnarray}  \label{test fd-nr}
   \zeta_0(t): =\begin{cases}
   0, & t\leq 0, \\[2mm]
 t^{\sqrt{1+4\mu}}, &   t\in[0,1),\\[2mm]
1\ \ \,  &  t\in[1,  +\infty).
\end{cases}
\end{eqnarray} 
 Then we define 
$$w_n(t):=\begin{cases}\nu_0 n^{\frac12-\tau_0}\zeta_0(t/n),& t \geq 1, \\[2mm]
\nu_0 n^{\frac12-\tau_0}\zeta_0(1/n)t,& t\leq 1,
\end{cases}
 \ \mbox{
with }\ \nu_0=(2^{2-2\tau_0}\pi)^{-1/2}(1-2\tau_0)^{-1/2}.$$
Let us note that $\mu\in (-1/4, 0)$ implies $2\tau_0\in (0,1)$. Therefore
  \begin{eqnarray*} 
  &&2^{2-2\tau_0}\pi \int^\infty_ {-t_\epsilon}  w_n'(t)^2(t+t_\epsilon)^{2\tau_0} dt\\
 &=&   (1-2\tau_0) n^{-1+2\tau_0}\Big[\int_\epsilon^1 (t+t_\epsilon)^{2\tau_0}   dt+ \int^n_1 t^{-4\tau_0} (t+t_\epsilon)^{2\tau_0}   dt\Big]
  \\[2mm]& \leq& (1-2\tau_0) n^{-1+2\tau_0}\Big[(1+t_\epsilon)^{2\tau_0}+\int^n_1 (t+t_\epsilon)^{-2\tau_0}   dt\Big] \\
  [2mm] &\leq&(1-2\tau_0) n^{-1+2\tau_0}(1+t_\epsilon)^{2\tau_0}+\big(1+t_\epsilon/n)^{1-2\tau_0}\\
  [2mm]&=&1+O(n^{-1+2\tau_0}).
    \end{eqnarray*}
   By (\ref{con V-l-c}) 
   \begin{eqnarray*} 
&& \Big| 2\pi\mu \int^1_0 \Big(V(r)r^2(-\ln r)^2-1\Big) \frac{(-\ln r)^{2\tau_0}}{r^2(-\ln r)^2} v_n(r)^2   rdr\Big| 
  \\[2mm] &\leq & 2\pi\mu \epsilon \nu_0^2 n^{1-2\tau_0} \Big[\int_\epsilon^1[\zeta_0(1/n)]^2 (t+t_\epsilon)^{-2+2\tau_0}dt+ \int^{\infty}_1  [\zeta_0(t/n)]^2 (t+t_\epsilon)^{-2+2\tau_0}dt\Big]
    \\[2mm] &\leq & 2\pi\mu  \epsilon \nu_0^2 n^{1-2\tau_0}   \Big( \int^n_0     n^{-2 \sqrt{1+4\mu}}  t^{-2\tau_0 }  dt + \int^{+\infty}_n     t^{-2+2\tau_0} dt\Big)
    \\[2mm] &= &  \frac{4\pi\mu  \nu_0^2}{1-2\tau_0} \epsilon.  
  \end{eqnarray*}
Thus 
$$\|u_n\|_{V,\mu}^2\leq 1+C(\epsilon+n^{-1+2\tau_0}).  $$ 
 \medskip
 
 Set $\hat u_n:=u_n/\|u_n\|_{V,\mu}$ and $\hat w_n=w_n/\|u_n\|_{V,\mu}$. Then, for all large $n$,
 \begin{eqnarray*} 
\tilde\cN(\hat w_n)
 &=&  \pi   \int^{\infty}_{-t_\epsilon-2\ln r_0}   e^{ \alpha  2^{-2\tau_0} (t+t_\epsilon)^{2\tau_0}\hat w_n^2-(t+t_\epsilon) }  dt
 \\[2mm]&>&  \pi  \int^{n}_{(1-\epsilon_0)n}    e^{\frac{ \alpha  2^{-2\tau_0}\nu_0^2n^{-1+2\tau_0}  t^{2-2\tau_0}}{1+C(\epsilon+n^{-1+2\tau_0})}-(t+t_\epsilon) }  dt
  \\[2mm]&\geq &  \pi e^{n\big[\frac{ \alpha  2^{-2\tau_0}\nu_0^2(1-\epsilon_0)^{2-2\tau_0}}{1+C(\epsilon+n^{-1+2\tau_0})}-(1+\frac {t_\epsilon}n) \big]}  \int^{n}_{(1-\epsilon_0)n} 1 dt  
  \\[2mm]&\geq & \pi (1-\epsilon_0)n
 \to   +\infty \quad{\rm as}\ \, n\to+\infty
    \end{eqnarray*}  
provided that
 $$
 \frac{ \alpha  2^{-2\tau_0}\nu_0^2(1-\epsilon_0)^{2-2\tau_0}}{1+C(\epsilon+n^{-1+2\tau_0})}-(1+\frac {t_\epsilon}n)
 \geq0,$$
   which holds for all large $n$ if
 $$
  \alpha >  2^{2\tau_0} \nu_0^{-2} = 4\pi  \sqrt{1+4\mu}  
$$ 
and  $\epsilon>0$, $\epsilon_0>0$ are chosen small enough.
 \hfill$\Box$\medskip

  \begin{lemma}\label{lm super 0}
Let   $\mu>0$,  $V:(0,1)\to[0,+\infty)$ be a continuous  function.
Then  for any $x_0\in B_1$ and $r_0\in (0,1)$ such that $\overline{B_{r_0}(x_0)} \subset B_1\setminus \{0\}$,
there exists a sequence of functions $\{u_n\}$ with support in $B_{r_0}(x_0)$ and radially symmetric about $x_0$, such that
 $\|u_n\|_{_V,\mu}= 1$ and
$$\lim_{n\to\infty}\int_{B_{1}} e^{\alpha |u_n|^2}dx=\infty \mbox{ for any } \alpha>4\pi.
$$
\end{lemma}
\noindent{\bf Proof.}  
Let $u_n$ be a sequence of nonnegative functions of $r=|x-x_0|$ with supports contained in $B_{r_0}(x_0)$, and extended by 0 outside the supporting sets. 
 Then 
\begin{eqnarray*} 
 \|u_n\|_{V,\mu}^2 \leq\cI(u_n): =  2\pi \int^{r_0}_0  u_n'(r)^2 rdr + 2\pi  \mu \|V\|_{L^\infty(B_{r_0}(x_0))}   \int^{r_0}_0  u_n^2(r)r  dr. 
    \end{eqnarray*}

 Let
  $$  r=e^{-\frac12 t}   \quad {\rm and}\quad u_n(r)=w_n(t),\quad t\in[-2\ln r_0,+\infty).$$ 
Then
  $$
\int_{B_{r_0}} e^{\alpha |u_n|^2}dx= 2\pi \int^{r_0}_0  e^{\alpha u_n^2}r dr = \pi   \int^\infty_{-2\ln r_0}  e^{ \alpha   w_n^2- t  }  dt:=\cN(w_n)$$
 and
 \begin{eqnarray*} 
    \cI(u_n) = 4\pi \int^\infty_ {-2\ln r_0}  w_n'(t)^2   dt+2\pi \mu \|V\|_{L^\infty(B_{r_0}(x_0))}   \int^\infty_ {-2\ln r_0}  w_n^2(t)e^{-t} dt.
    \end{eqnarray*}

Next we  construct the functions $\{w_n\}$.   Let
 \begin{eqnarray}  \label{test fd-nr}
   \zeta_0(t): =\left\{ \arraycolsep=1pt
\begin{array}{lll}
0\ \ \, &{\rm for}\ \    t\in[0, 1/4),\\[2mm]
2t-\frac12\ \ \, &{\rm for}\ \    t\in[1/4, 1/2),\\[2mm]
 t\ \ \, &{\rm for}\ \    t\in[1/2, 1),\\[2mm]
1\ \ \,  &{\rm for}\ \    t\in[1,  +\infty).
\end{array}
\right.  
\end{eqnarray} 
We then define, for all large $n$,
$$w_n(t):= \sqrt{\frac n{4\pi}}\, \zeta_0\big(\frac t n\big). $$
Clearly, for all large $n$,
  \begin{eqnarray*} 
   4\pi  \int^\infty_ {-2\ln r_0}   w_n'(t)^2  dt=      n^{-1}\Big(\int^{n/4}_{n/2} 2dt+\int_{n/2}^n 1    dt\Big)=1
    \end{eqnarray*}
    and 
   \begin{eqnarray*} 
  2\pi \int^\infty_{-2\ln r_0}   w^2_n(t) e^{-t} dt &=&      2\pi    \int^\infty_{n/4}   w_n^2(t) e^{-t} dt      \\[2mm]
  &  \leq&     \frac 2{n}  \int^{n}_{n/4}   t^2 e^{-t} dt + \frac n2\int^{+\infty}_ {n}    e^{-t} dt
\\[2mm]&   \leq & e^{-n/2}.
    \end{eqnarray*}
Thus 
$$\|u_n\|_{V,\mu}^2\leq \cI_n(u_n)\leq 1+   2\pi \mu \|V\|_{L^\infty(B_{r_0}(x_0))}  e^{-n/2}=:1+Ce^{-n/2}.  $$ 
 \medskip
 
 Set $\hat u_n:=u_n/\|u_n\|_{V,\mu}$ and $\hat w_n=w_n/\|u_n\|_{V,\mu}$. Then, for all large $n$ and $\epsilon_0\in (0, 1/2)$,
 \begin{eqnarray*} 
\cN(\hat w_n)
 &=&  \pi   \int^{+\infty}_0   e^{ \alpha   \hat w_n^2(t)- t  }  dt
 \\[2mm]&>&\pi    \int^{n}_{(1-\epsilon_0)n}   e^{ \frac{\alpha }{1+Ce^{-n/2}}\frac{n}{4\pi} (\frac tn )^2- t   }  dt
    \\[2mm]&\geq &  \pi  e^{n\big[\frac{\alpha (1-\epsilon_0)^2}{4\pi(1+Ce^{-n/2})} -1\big]  }   \int^{n}_{(1-\epsilon_0)n} 1 dt  \,\geq  \pi     \epsilon_0 n \to \infty \quad{\rm as}\ \, n\to \infty
    \end{eqnarray*}  
provided that $$
\frac{\alpha (1-\epsilon_0)^2}{4\pi(1+Ce^{-n/2})} -1
\geq 0,$$ 
which holds for all large $n$ if
$\alpha>4\pi \mbox{ and $ \epsilon_0>0$ is chosen sufficiently small}. $
    \hfill$\Box$ \medskip

\noindent{\bf Proof of Theorem \ref{cr 0}. } {\it Part $(i)$-$(ii)$: } In the radially symmetric 
case, the bound  follows from  Lemma \ref{teo 3.1}  and in the nonradial case, the bound  follows from   \cite[Theorem 1]{M}.  

 Part $(iii)$-$(iv)$ follows from Lemma \ref{lm super c}   for $\alpha> m_\mu$ in the radial case and from Lemma \ref{lm super 0} for $\alpha> 4\pi$ in the nonradial case. 
 \hfill$\Box$\bigskip

\noindent{\bf Proof of Theorem \ref{cr 1}. }  For   $u\in \cH^1_{_V,\mu, 0}(B_1)$ such that 
$\|u\|_{_V,\mu}\leq 1$, using (\ref{con V}) and $\mu<0$,
   we obtain that  $u\in \cH^1_{\mu, 0}(B_1)$  and $\|u\|_{\mu}\leq \|u\|_{_V,\mu}\leq 1$.

{\it Part $(i)$. }  If $u$ is radially symmetric, it follows from
Lemma \ref{teo 3.1} that 
$$\int_{B_1} e^{m_\mu |u|^2}dx<\infty .$$

 {\it Part $(ii)$. }  Without loss of generality we assume    $u\in C_c^{0,1}(B_1)$,
 and  let $u^*$ denote its associated radially decreasing rearrangement.  

By the P\'olya-Szeg\H{o} inequality and rearrangement inequalities we have, since  $r\in(0,1)\to V(r)$ is decreasing,
$$\int_{B_1} |\nabla u^*|^2 dx\leq \int_{B_1} |\nabla u |^2 dx,$$  
\begin{eqnarray*}  
\int_{B_1}   u^2 V(|x|)  dx    \leq  \int_{B_1}   (u^*)^2 V^*(|x|)  dx = \int_{B_1} (u ^*)^2 V(|x|) dx 
   \end{eqnarray*}
   and 
     \begin{eqnarray*}    
\int_{B_1}  e^{\alpha u^p}  dx \leq  \int_{B_1} e^{\alpha (u^*)^p}  dx. 
  \end{eqnarray*}
  Hence, for $\mu\in[ -\frac14,0)$, we have 
\begin{eqnarray*} 
  1  \geq    \int_{B_1} \Big(|\nabla u|^2+\mu u(x)^2 V(|x|) \Big)   dx
 \geq        \int_{B_1} \Big(|\nabla u^*|^2+\mu   u^*(x)^2V(|x|) \Big)   dx.
    \end{eqnarray*} 
  Now we can apply Lemma \ref{teo 3.1} to $u^*$  to obtain the desired conclusion.
   \smallskip

  Part $(iii)$ follows from Lemma \ref{lm super c}.  \hfill$\Box$ \medskip

 \subsection{The case $\mu=-\frac14$} We give the proof of
 Theorem \ref{cr 2} here.
 
   For any  function $u\in \cH^1_{_V, -\frac14, 0}(B_1)$ such that 
$\|u\|_{_V,-\frac14}\leq 1$, by (\ref{con V}),
     $$u\in \cH_{0}(B_1)\quad{\rm  and}\quad \|u\|_{-\frac 14}\leq \|u\|_{_V,-\frac14}\leq 1.$$

{\it Part $(i)$. }  If $u$ is radial, it follows from  Lemma \ref{teo 3.2}
that for any $\alpha>0$, $p\in(0,1)$, there exists $C=C_{\alpha, p}>0$ independent of $u$ such that 
$$\int_{B_1} e^{\alpha |u|^p}dx\leq C .$$

{\it Part $(ii)$. }      As before we may assume  $u\in C_c^{1}(B_1)$,
 and  let $u^*$ denote its associated radially decreasing rearrangement.   By the arguments in
 the proof of Theorem \ref{cr 1}$(ii)$,  we have 
\begin{eqnarray*} 
     \int_{B_1} \Big(|\nabla u|^2 -\frac14 u(x)^2 V(|x|) \Big)   dx
 \geq        \int_{B_1} \Big(|\nabla u^*|^2-\frac14   u^*(x)^2V(|x|) \Big)   dx.
    \end{eqnarray*} 
  Now we can apply Lemma \ref{teo 3.2} to $u^*$  to obtain the desired bound.  \smallskip

{\it Part $(iii)$. }  
By (\ref{con V-l2}), there exist $C>0$ and $\theta>0$ such that
\[
|V(r) {r^2(-\ln r)^2}-1|\leq C(-\ln r)^{-\theta} \quad \mbox{for}\ r\in(0,\frac14).
\]

Let 
 $$ \begin{cases} u (r)=(-\ln r)^{\frac12} v (r)\quad{\rm for}\ \, r\in(0,1),\\[2mm]
  \frac 14=e^{\frac12(1-  e^{t_0}) },\  r=e^{\frac12(1-  e^{t+t_0}) } ,\  v (r)=w (t);
  \end{cases}$$
 then 
 $$ dr=-\frac12e^{\frac12(1-  e^{t+t_0}) } e^{t+t_0}dt,\ \lim_{t\to-t_0^+} w (t)=\lim_{r\to1^-}v (r)=0  $$
and  
 \begin{eqnarray*} 
 \|u\|_{-1/4}^2  =    2\pi \int^1_0   v'(r)^2 r(-\ln r)  dr
 =  2\pi \int^{+\infty}_{-t_0}  w'(t)^2 \,  \frac{e^{t+t_0}-1}{e^{t+t_0}}   dt, 
    \end{eqnarray*}
     \begin{eqnarray*} 
  \Big|   \int^{1/4}_0 \Big(V(r)r^2(-\ln r)^2-1\Big) \frac{-\ln r}{r^2(-\ln r)^2} v(r)^2   rdr\Big| 
  &\leq&  C   \int^{1/4}_0  \frac{1}{r(-\ln r)^{1+\theta}} v(r)^2   dr 
      \\[2mm] &= &  C   \int^{+\infty}_{0} 
      w(t)^2 \,  \frac{e^{t+t_0}}{(e^{t+t_0}-1)^{1+\theta}}  dt ,   
  \end{eqnarray*}   
  \begin{eqnarray*}  
  \int_{B_1} e^{\alpha |u|^p}dx= 2\pi \int^1_0  e^{\alpha |u|^p } r dr
 =   \pi   \int^{+\infty}_{-t_0}  e^{\alpha \frac{(e^{t+t_0}-1)^p}{2^p} |w|^p+1-e^{t+t_0} +(t+t_0)}  dt.
   \end{eqnarray*}

We now construct a sequence of radially symmetric functions  $\{\hat u_n\}\subset C_c^\infty(B_1)$, whose supporting sets shrink to the origin as $n\to\infty$, and satisfies 
\[
\|\hat u_n\|_{V,-1/4}=1,
 \ \lim_{n\to\infty}\int_{B_1} e^{\alpha |\hat u_n|^p}dx=\infty
 \]
for any $p\geq 1$ and any $\alpha>0$.  

 Given $n>4$, we let $\eta_n:(-\infty, +\infty)\to[0,1]$ be a  smooth function such that 
\begin{eqnarray} \label{cutoff f1-} 
  \eta_n(t) =\left\{ \arraycolsep=1pt
\begin{array}{lll}
 0\ \ \ \ &{\rm for}\ \,   t\in(-\infty,0),
\\[2mm]
1\ \  &{\rm for}\ \,   t\in(1,n-1),\\[2mm]
0\ \  &{\rm for}\ \,   t\in(n,+\infty).
\end{array}
\right.  
    \end{eqnarray} 
Define
   $$b_1(n):=\int^{n}_1  \eta_n (s) ds\quad{\rm and}\quad b_2(n):=\sqrt{4\pi \int^{n}_0  \eta_n(s)^2 ds};$$
then 
\begin{eqnarray*}  
 \frac{b_1(n)}{b_2(n)}  \geq     \frac{ n-2  }{\sqrt{4\pi n} } \to+\infty\quad{\rm as}\ \, n\to+\infty.
    \end{eqnarray*} 
Set
\begin{eqnarray} \label{cutoff f2-c}  
\phi_n(t):=  \eta_n(t)-\eta_n( t-2n) \quad{\rm for}\ \, t\in \R,
  \end{eqnarray} 
  and
  $$w_n(t):=\frac1{b_2(n)+1} \int_0^t\phi_n (s)ds\quad{\rm for}\ \, t\geq0; $$
   then $w_n$ is nonnegative and $w_n\in C_c^\infty((0,\infty))$.

 Moreover, 
 \begin{eqnarray*}  
\|u_n\|_{-1/4}^2= 2\pi \int^{+\infty}_{-t_0}  w_n '(t)^2 \,  \frac{e^{t+t_0}-1}{e^{t+t_0}}   dt &\leq &    \frac{2\pi}{[b_2(n)+1]^2}   \int^{+\infty}_0  \phi^2_n(t)  dt 
\\[2mm]&  =&   \frac{b_2(n)^2 }{[b_2(n)+1]^2}<1,
  \end{eqnarray*}
  and due to $w_n(t)=0$ for $t\leq 0$, 
   \begin{eqnarray*} 
  \|u_n\|_{V, -1/4}^2-\|u_n\|_{-1/4}^2&=&   \int^1_0 \Big(V(r)r^2(-\ln r)^2-1\Big) \frac{-\ln r}{r^2(-\ln r)^2} v_n(r)^2   rdr \\
  &=&    \int^{1/4}_0 \Big(V(r)r^2(-\ln r)^2-1\Big) \frac{1}{r(-\ln r)} v_n(r)^2   dr \\
  &\leq&  C \int^{+\infty}_{0}  w_n(t)^2 \,  \frac{e^{t+t_0}}{(e^{t+t_0}-1)^{1+\theta}}   dt 
      \\[2mm] &\leq & C \int_{0}^{+\infty}    \frac{t^2e^{t+t_0}}{ (e^{t+t_0}-1)^{1+\theta} }    dt
      \\[2mm] &=:&C_0. 
  \end{eqnarray*}   
 Therefore
\[
\|u_n\|_{V, -1/4}^2\leq 1+C_0.
\]
  
 Note that for $ t\in[n, 2n]$,
 $$w_n(t)=\int_0^{n} \phi_n (s)ds=\frac{b_1(n)}{b_2(n)+1}$$
   and with
   \[\hat u_n:=u_n/\|u_n\|_{V,-1/4},\ \hat w_n:=w_n/\|u_n\|_{V, -1/4},
   \]
   we have, for all large $n$,
 \begin{eqnarray*} 
 \int_{B_1} e^{\alpha |\hat u_n|^p}dx&=&
   \pi   \int^{+\infty}_{-t_0} e^{\alpha \big(\frac{ e^{t+t_0}-1 }{2 }\big)^p |\hat w_n |^p+1-e^{t+t_0} +(t+t_0) }  dt\\
   &>& \pi   \int^{2n}_{n}   e^{\alpha \big(\frac{ e^{t+t_0}-1 }{2 }\big)^p |\hat w_n|^p+1-e^{t+t_0} +(t+t_0) }  dt
\\[2mm]  &\geq &  \pi  \int^{2n}_{n}     e^{\alpha \big[\frac{b_1(n)}{b_2(n)+1} \frac{ 1 }{4(1+C_0)^{2/p}} \big]^pe^{p(t+t_0)}-e^{t+t_0}  }  dt
 \\[2mm]&\to & +\infty\quad{\rm as}\ \ n\to+\infty, 
     \end{eqnarray*}   
     since $\frac{b_1(n)}{b_2(n)+1}\to\infty$ as $n\to\infty$.
 This completes the proof.  \hfill$\Box$

\end{document}